\documentclass[]{article}

\usepackage[a4paper,top=2.1cm,bottom=2.60cm,left=2.1cm,right=2.1cm]{geometry} 

\usepackage{latexsym}
\usepackage{amsmath}
\usepackage{amsfonts}
\usepackage{mathrsfs}
\usepackage{amstext}
\usepackage{amssymb}
\usepackage{amsthm}       

\usepackage{moreverb}

\usepackage[dvips,colorlinks,bookmarksopen,bookmarksnumbered,citecolor=red,urlcolor=red]{hyperref}

\newcommand\BibTeX{{\rmfamily B\kern-.05em \textsc{i\kern-.025em b}\kern-.08em
T\kern-.1667em\lower.7ex\hbox{E}\kern-.125emX}}

\numberwithin{equation}{section}

\newtheorem{thm}[equation]{Theorem}
\newtheorem{prop}[equation]{Proposition}

\newtheorem{rem}[equation]{Remark}
\newtheorem{lem}[equation]{Lemma}

\begin{document}

\title{A singularly perturbed Dirichlet problem for the Laplace operator in a periodically perforated domain. A functional analytic approach}

\author{Paolo Musolino}

\date{}

\maketitle

\noindent
{\bf Abstract:}
Let $\Omega$ be a sufficiently regular bounded open connected subset of $\mathbb{R}^n$ such that $0 \in \Omega$ and that $\mathbb{R}^n \setminus \mathrm{cl}\Omega$ is connected. Then we take $q_{11},\dots, q_{nn}\in  ]0,+\infty[$ and $p \in Q\equiv \prod_{j=1}^{n}]0,q_{jj}[$. If $\epsilon$ is a small positive number, then we define the periodically perforated domain $\mathbb{S}[\Omega_\epsilon]^{-} \equiv \mathbb{R}^n\setminus \cup_{z \in \mathbb{Z}^n}\mathrm{cl}\bigl(p+\epsilon \Omega +\sum_{j=1}^n (q_{jj}z_j)e_j\bigr)$, where $\{e_1,\dots,e_n\}$ is the canonical basis of $\mathbb{R}^n$. For $\epsilon$ small and positive, we introduce a particular Dirichlet problem for the Laplace operator in the set $\mathbb{S}[\Omega_\epsilon]^{-}$. Namely, we consider a Dirichlet condition on the boundary of the set $p+\epsilon \Omega$, together with a periodicity condition. Then we show real analytic continuation properties of the solution and of the corresponding energy integral as functionals of the pair of $\epsilon$ and of the Dirichlet datum on $p+\epsilon \partial \Omega$, around a degenerate pair with $\epsilon=0$. 

\vspace{11pt}

\noindent
{\bf MOS:} 35 J 25; 31 B 10; 45 A 05; 47 H 30


\maketitle

\section{Introduction}
In this article, we consider a Dirichlet problem in a periodically perforated domain with small holes. We fix once for all  a natural number
\[
n\in {\mathbb{N}}\setminus\{0,1 \}
\]
and
\[
 (q_{11},\dots,q_{nn})\in]0,+\infty[^{n}
\]
and a periodicity cell
\[
Q\equiv\Pi_{j=1}^{n}]0,q_{jj}[\,.
\]
Then we denote by $q$ the diagonal matrix
\[
q\equiv \left(
\begin{array}{cccc}
q_{11} &   0 & \dots & 0   
\\
0          &q_{22} &\dots & 0
\\
\dots & \dots & \dots & \dots  
\\
0& 0 & \dots & q_{nn}
\end{array}\right),
\]
by $|Q|$ the measure of the fundamental cell $Q$, and by $\nu_Q$ the outward unit normal to $\partial Q$, where it exists. Clearly, 
\[
q {\mathbb{Z}}^{n}\equiv  \{qz:\,z\in{\mathbb{Z}}^{n}\}
\]
is the set of vertices of a periodic subdivision of ${\mathbb{R}}^{n}$ corresponding to the fundamental cell $Q$. Let
\[
m\in {\mathbb{N}}\setminus\{0\}\,,\qquad\alpha\in]0,1[\,.
\]
Then we take a point $p \in Q$ and a bounded open connected subset $\Omega$ of ${\mathbb{R}}^{n}$ of class $C^{m,\alpha}$ such that $\Omega^{-}\equiv {\mathbb{R}}^{n}\setminus{\mathrm{cl}}\Omega$ is connected and that $0 \in \Omega$. If $\epsilon \in \mathbb{R}$, then we set
\[
\Omega_\epsilon\equiv p +\epsilon \Omega\,.
\]
Then we take $\epsilon_0>0$ such that $\mathrm{cl}\Omega_\epsilon \subseteq Q$ for $|\epsilon|< \epsilon_0$, and we introduce the periodically perforated domain
\[
{\mathbb{S}} [\Omega_\epsilon]^{-}\equiv {\mathbb{R}}^{n}\setminus\cup_{z \in \mathbb{Z}^n}{\mathrm{cl}}(\Omega_\epsilon+qz)\,,
\]
for $\epsilon \in ]-\epsilon_0,\epsilon_0[$. Next we fix a function $g_0 \in C^{m,\alpha}(\partial \Omega)$. For each pair $(\epsilon,g) \in ]0,\epsilon_0[\times C^{m,\alpha}(\partial \Omega)$ we consider the Dirichlet problem 
 \begin{equation}\label{bvp:Direps}
 \left \lbrace 
 \begin{array}{ll}
 \Delta u (x)= 0 & \textrm{$\forall x \in {\mathbb{S}} [\Omega_\epsilon]^{-}$}\,, \\
u(x+qe_i) =u(x) &  \textrm{$\forall x \in \mathrm{cl} {\mathbb{S}} [\Omega_\epsilon]^{-}\,, \quad \forall i \in \{1,\dots,n\}$}\,, \\
u(x)=g\bigl(\frac{1}{\epsilon}(x-p)\bigr) & \textrm{$\forall x \in \partial \Omega_\epsilon$}\,,
 \end{array}
 \right.
 \end{equation}
where $\{e_1,\dots,e_n\}$ is the canonical basis of $\mathbb{R}^n$. If $(\epsilon,g) \in ]0,\epsilon_0[\times C^{m,\alpha}(\partial \Omega)$, then problem \eqref{bvp:Direps} has a unique solution in $C^{m,\alpha}(\mathrm{cl} {\mathbb{S}} [\Omega_\epsilon]^{-})$, and we denote it by $u[\epsilon,g](\cdot)$ (cf. Proposition \ref{prop:Dirsol}.)

Then we pose the following questions:
\begin{enumerate}
\item[(i)] Let $x$ be fixed in $\mathbb{R}^n \setminus (p+q\mathbb{Z}^n)$. What can be said on the map $(\epsilon,g) \mapsto u[\epsilon,g](x)$ around $(\epsilon,g)=(0,g_0)$?
\item[(ii)] What can be said on the map $(\epsilon,g) \mapsto\int_{Q\setminus \mathrm{cl} \Omega_\epsilon}|D_x u[\epsilon,g](x)|^2 \, dx$ around $(\epsilon,g)=(0,g_0)$?
\end{enumerate}

Questions of this type have long been investigated, \textit{e.g.}, for problems on a bounded domain with a small hole with the methods of asymptotic analysis, which aims at giving complete asymptotic expansions of the solutions in terms of the parameter $\epsilon$. It is perhaps difficult to provide a complete list of the contributions. Here, we mention the work of Ammari and Kang \cite[Ch.~5]{AmKa07}, Ammari, Kang, and Lee \cite[Ch.~3]{AmKaLe09}, Kozlov, Maz'ya, and Movchan \cite{KoMaMo99}, Maz'ya, Nazarov, and Plamenewskij \cite{MaNaPl00i, MaNaPl00ii}, Ozawa \cite{Oz83}, Vogelius and Volkov \cite{VoVo00}, Ward and Keller \cite{WaKe93}. We also mention the vast literature of homogenization theory (cf. \textit{e.g.}, Dal Maso and Murat \cite{DaMu04}.)

Here instead we wish to characterize the behaviour of $u[\epsilon,g](\cdot)$ at $(\epsilon,g)=(0,g_0)$ by a different approach. Thus for example, if we consider a certain functional, say $f(\epsilon,g)$, relative to the solution such as, for example, one of those considered in questions (i)-(ii) above, we would try to prove that $f(\cdot,\cdot)$ can be continued real analytically around $(\epsilon,g)=(0,g_0)$. We observe that our approach does have certain advantages (cf.~\textit{e.g.}, Lanza \cite{La07a}.) Such a project has been carried out by Lanza de Cristoforis in several papers for problems in a bounded domain with a small hole (cf.~\textit{e.g.}, Lanza \cite{La02, La04,La08,La10}.) In the frame of linearized elastostatics, we also mention, \textit{e.g.}, Dalla Riva and Lanza \cite{DaLa10, DaLa10b}.

As far as problems in periodically perforated domains are concerned, we mention, for instance, the work of Ammari, Kang, and Touibi \cite{AmKaTo05}, where a linear transmission problem is considered in order to compute an asymptotic expansion of the  effective electrical conductivity of a periodic dilute composite (see also Ammari and Kang \cite[Ch. 8]{AmKa07}.) Furthermore, we note that periodically perforated domains are extensively studied in the frame of homogenization theory. Among the vast literature, here we mention, \textit{e.g.}, Cioranescu and Murat \cite{CiMu82i, CiMu82ii}, Ansini and Braides \cite{AnBr02}. We also observe that boundary value problems in domains with periodic inclusions can be analysed, at least for the two dimensional case, with the method of functional equations. Here we mention, \textit{e.g.}, Mityushev and Adler \cite{MiAd02i}, Rogosin, Dubatovskaya, and Pesetskaya \cite{RoDuPe09}, Castro and Pesetskaya \cite{CaPe10}.

We now briefly outline our strategy. We first convert problem \eqref{bvp:Direps} into an integral equation by exploiting potential theory. Then we observe that the corresponding integral equation can be written, after an appropriate rescaling, in a form which can be analysed by means of the Implicit Function Theorem around the degenerate case in which $(\epsilon,g)=(0,g_0)$, and we represent the unknowns of the integral equation in terms of $\epsilon$ and $g$. Next we exploit the integral representation of the solutions, and we deduce the representation of $u[\epsilon,g](\cdot)$ in terms of $\epsilon$ and $g$.

This article is organized as follows. Section \ref{not} is a section of preliminaries. In Section \ref{form}, we formulate problem \eqref{bvp:Direps} in terms of an integral equation and we show that the solutions of the integral equation depend real analytically on $\epsilon$ and $g$. In Section \ref{rep}, we show that the results of Section \ref{form} can be exploited to prove our main Theorem \ref{thm:rep} on the representation of $u[\epsilon,g](\cdot)$, and Theorem \ref{thm:en} on the representation of the energy integral of $u[\epsilon,g](\cdot)$ on a perforated cell. At the end of this article, we have enclosed an Appendix with some results exploited in the paper.

\section{Preliminaries and notation}\label{not}

We now introduce the notation in accordance with Lanza \cite[p.~66]{La08}.

We  denote the norm on 
a   normed space ${\mathcal X}$ by $\|\cdot\|_{{\mathcal X}}$. Let 
${\mathcal X}$ and ${\mathcal Y}$ be normed spaces. We endow the  
space ${\mathcal X}\times {\mathcal Y}$ with the norm defined by 
$\|(x,y)\|_{{\mathcal X}\times {\mathcal Y}}\equiv \|x\|_{{\mathcal X}}+
\|y\|_{{\mathcal Y}}$ for all $(x,y)\in  {\mathcal X}\times {\mathcal 
Y}$, while we use the Euclidean norm for ${\mathbb{R}}^{n}$.
 For 
standard definitions of Calculus in normed spaces, we refer to 
Prodi and Ambrosetti~\cite{PrAm73}. The symbol ${\mathbb{N}}$ denotes the 
set of natural numbers including $0$.  The inverse function of an 
invertible function $f$ is denoted $f^{(-1)}$, as opposed to the 
reciprocal of a real-valued function $g$, or the inverse of a 
matrix $A$, which are denoted $g^{-1}$ and $A^{-1}$, respectively. A 
dot ``$\cdot$'' denotes the inner product in ${\mathbb R}^{n}$. Let $A$ be a 
matrix. Then $A^{t}$ denotes the transpose matrix of $A$  and 
 $A_{ij}$ denotes 
the $(i,j)$-entry of $A$. If $A$ is invertible, 
we set $A^{-t}\equiv \left(A^{-1}\right)^{t}$. 
 Let 
${\mathbb{D}}\subseteq {\mathbb {R}}^{n}$. Then $\mathrm{cl}\,{\mathbb{D}}$ 
denotes the 
closure of ${\mathbb{D}}$ and $\partial{\mathbb{D}}$ denotes the boundary of ${\mathbb{D}}$. For all $R>0$, $ x\in{\mathbb{R}}^{n}$, 
$x_{j}$ denotes the $j$-th coordinate of $x$, 
$| x|$ denotes the Euclidean modulus of $ x$ in
${\mathbb{R}}^{n}$, and ${\mathbb{B}}_{n}( x,R)$ denotes the ball $\{
y\in{\mathbb{R}}^{n}:\, | x- y|<R\}$. The symbol $\mathrm{id}_n$ denote the identity map from $\mathbb{R}^n$, \textit{i.e}, $\mathrm{id}_n(x)=x$ for all $x \in \mathbb{R}^n$. If $z \in \mathbb{C}$, then $\overline{z}$ denotes the conjugate complex number of $z$. 
Let $\Omega$ be an open 
subset of ${\mathbb{R}}^{n}$. The space of $m$ times continuously 
differentiable real-valued functions on $\Omega$ is denoted by 
$C^{m}(\Omega,{\mathbb{R}})$, or more simply by $C^{m}(\Omega)$. 
${\mathcal{D}}(\Omega)$ denotes the space of functions of  $C^{\infty}(\Omega)$
with compact support. The dual ${\mathcal{D}}'(\Omega)$ denotes the space
of distributions in $\Omega$.
Let $r\in {\mathbb{N}}\setminus\{0\}$.
Let $f\in \left(C^{m}(\Omega)\right)^{r}$. The 
$s$-th component of $f$ is denoted $f_{s}$, and $Df$ denotes the Jacobian matrix
$\left(\frac{\partial f_s}{\partial
x_l}\right)_{  \substack{ s=1,\dots ,r,    \\  l=1,\dots ,n}       }$.
Let  $\eta\equiv
(\eta_{1},\dots ,\eta_{n})\in{\mathbb{N}}^{n}$, $|\eta |\equiv
\eta_{1}+\dots +\eta_{n}  $. Then $D^{\eta} f$ denotes
$\frac{\partial^{|\eta|}f}{\partial
x_{1}^{\eta_{1}}\dots\partial x_{n}^{\eta_{n}}}$.    The
subspace of $C^{m}(\Omega )$ of those functions $f$ whose derivatives $D^{\eta }f$ of
order $|\eta |\leq m$ can be extended with continuity to 
$\mathrm{cl}\,\Omega$  is  denoted $C^{m}(
\mathrm{cl}\,\Omega )$. 
The
subspace of $C^{m}(\mathrm{cl}\,\Omega ) $  whose
functions have $m$-th order derivatives that are
H\"{o}lder continuous  with exponent  $\alpha\in
]0,1]$ is denoted $C^{m,\alpha} (\mathrm{cl}\,\Omega )$  
(cf.~\textit{e.g.},~Gilbarg and 
Trudinger~\cite{GiTr83}.) The subspace of $C^{m}(\mathrm{cl}\,\Omega ) $ of those functions $f$ such that $f_{|{\mathrm{cl}}(\Omega\cap{\mathbb{B}}_{n}(0,R))}\in
C^{m,\alpha}({\mathrm{cl}}(\Omega\cap{\mathbb{B}}_{n}(0,R)))$ for all $R\in]0,+\infty[$ is denoted $C^{m,\alpha}_{{\mathrm{loc}}}(\mathrm{cl}\,\Omega ) $.  Let 
${\mathbb{D}}\subseteq {\mathbb{R}}^{r}$. Then $C^{m
,\alpha }(\mathrm{cl}\,\Omega ,{\mathbb{D}})$ denotes
$\left\{f\in \left( C^{m,\alpha} (\mathrm{cl}\,\Omega )\right)^{r}:\ f(
\mathrm{cl}\,\Omega )\subseteq {\mathbb{D}}\right\}$. \par
Now let $\Omega $ be a bounded
open subset of  ${\mathbb{R}}^{n}$. Then $C^{m}(\mathrm{cl}\,\Omega )$ 
and $C^{m,\alpha }({\mathrm{cl}}\,
\Omega )$ are endowed with their usual norm and are well known to be 
Banach spaces  (cf.~\textit{e.g.}, Troianiello~\cite[\S 1.2.1]{Tr87}.) 
We say that a bounded open subset $\Omega$ of ${\mathbb{R}}^{n}$ is of class 
$C^{m}$ or of class $C^{m,\alpha}$, if it is a 
manifold with boundary imbedded in 
${\mathbb{R}}^{n}$ of class $C^{m}$ or $C^{m,\alpha}$, respectively
 (cf.~\textit{e.g.}, Gilbarg and Trudinger~\cite[\S 6.2]{GiTr83}.) 
We denote by 
$
\nu_{\Omega}
$
the outward unit normal to $\partial\Omega$.  For standard properties of functions 
in Schauder spaces, we refer the reader to Gilbarg and 
Trudinger~\cite{GiTr83} and to Troianiello~\cite{Tr87}
(see also  Lanza \cite[\S 2, Lem.~3.1, 4.26, Thm.~4.28]{La91}, 
 Lanza and Rossi \cite[\S 2]{LaRo04}.) \par
 
 We retain the standard notation of $L^p$ spaces and of corresponding norms.

\par
If $\mathbb{M}$ is a manifold  imbedded in 
${\mathbb{R}}^{n}$ of class $C^{m,\alpha}$, with $m\geq 1$, 
$\alpha\in ]0,1[$, one can define the Schauder spaces also on $\mathbb{M}$ by 
exploiting the local parametrizations. In particular, one can consider 
the spaces $C^{k,\alpha}(\partial\Omega)$ on $\partial\Omega$ for 
$0\leq k\leq m$ with $\Omega$ a bounded open set of class $C^{m,\alpha}$,
and the trace operator from $C^{k,\alpha}({\mathrm{cl}}\Omega)$ to
$C^{k,\alpha}(\partial\Omega)$ is linear and continuous. Moreover, for 
each $R>0$ such that ${\mathrm{cl}}\Omega\subseteq 
{\mathbb{B}}_{n}(0,R)$, there exists a linear and continuous 
extension operator from $C^{k,\alpha}(\partial\Omega)$ to 
$C^{k,\alpha}({\mathrm{cl}}\Omega)$, and of $C^{k,\alpha}({\mathrm{cl}}\Omega)$
to $C^{k,\alpha}({\mathrm{cl}}{\mathbb{B}}_{n}(0,R))$ (cf.~\textit{e.g.}, Troia\-niello~\cite[Thm.~1.3, Lem.~1.5]{Tr87}.)  We denote by $d\sigma$ the area element of a manifold imbedded in ${\mathbb{R}}^{n}$. \par

We note that 
throughout the paper ``analytic'' means ``real analytic''. For the 
definition and properties of analytic operators, we refer to Prodi and 
Ambrosetti~\cite[p. 89]{PrAm73} and to Deimling \cite[p.~150]{De85}. Here we just recall that if $\mathcal{X}$, $\mathcal{Y}$ are (real) Banach spaces, and if $F$ is an operator from an open subset $\mathcal{W}$ of $\mathcal{X}$ to $\mathcal{Y}$, then $F$ is real analytic in $\mathcal{W}$ if for every $x_0 \in \mathcal{W}$ there exist $r>0$ and continuous symmetric $n$-linear operators $A_n$ from $\mathcal{X}^n$ to $\mathcal{Y}$ such that $\sum_{n\geq 1} \|A_n\| r^n <\infty$ and $F(x_0+h)=F(x_0)+\sum_{n \geq 1} A_n(h,\dots,h)$ for $\|h\|_{\mathcal{X}} \leq r$ (cf. \textit{e.g.}, Prodi and 
Ambrosetti~\cite[p. 89]{PrAm73} and Deimling \cite[p.~150]{De85}.) In particular, we mention that the 
pointwise product in Schauder spaces is bilinear and continuous, and 
thus analytic, and that the map which takes a nonzero function 
to its reciprocal, or an invertible matrix of functions to its 
inverse matrix is real analytic in Schauder spaces 
(cf.~\textit{e.g.}, Lanza and Rossi \cite[pp.~141, 142]{LaRo04}.)\par

We denote by $S_{n}$  the function from  
${\mathbb{R}}^{n}\setminus\{0\}$ to
${\mathbb{R}}$ defined by 
\[
S_{n}(x)\equiv
\left\{
\begin{array}{lll}
\frac{1}{s_{n}}\log |x| \qquad &   \forall x\in 
{\mathbb{R}}^{n}\setminus\{0\},\quad & {\mathrm{if}}\ n=2\,,
\\
\frac{1}{(2-n)s_{n}}|x|^{2-n}\qquad &   \forall x\in 
{\mathbb{R}}^{n}\setminus\{0\},\quad & {\mathrm{if}}\ n>2\,,
\end{array}
\right.
\]
where $s_{n}$ denotes the $(n-1)$-dimensional measure of 
$\partial{\mathbb{B}}_{n}$. $S_{n}$ is well-known to be the 
fundamental solution of the Laplace operator.

If $y\in{\mathbb{R}}^{n}$ and $f$ is a function defined in ${\mathbb{R}}^{n}$, we set
$\tau_{y}f(x)\equiv f(x-y)$ for all $x\in {\mathbb{R}}^{n}$. If $u$ is a distribution in ${\mathbb{R}}^{n}$, then we set
\[
<\tau_{y}u,f>=<u,\tau_{-y}f>\qquad\forall f\in {\mathcal{D}}( {\mathbb{R}}^{n})\,.
\]
We denote by $E_{2\pi i q^{-1} z}$, the function defined by
\[
E_{2\pi i q^{-1} z}(x)\equiv e^{2\pi i (q^{-1} z)\cdot x}
\qquad  \forall x\in{\mathbb{R}}^{n}\,,
\]
for all $z\in {\mathbb{Z}}^{n}$.

If $\Omega$ is an open subset of ${\mathbb{R}}^{n}$, $k\in {\mathbb{N}}$, $\beta\in]0,1]$, we set
\[
C^{k}_{b}({\mathrm{cl}}\Omega)\equiv
\{
u\in C^{k}({\mathrm{cl}}\Omega):\,
D^{\gamma}u\ {\mathrm{is\ bounded}}\ \forall\gamma\in {\mathbb{N}}^{n}\
{\mathrm{such\ that}}\ |\gamma|\leq k
\}\,,
\]
and we endow $C^{k}_{b}({\mathrm{cl}}\Omega)$ with its usual  norm
\[
\|u\|_{ C^{k}_{b}({\mathrm{cl}}\Omega) }\equiv
\sum_{|\gamma|\leq k}\sup_{x\in {\mathrm{cl}}\Omega }|D^{\gamma}u(x)|\qquad\forall u\in C^{k}_{b}({\mathrm{cl}}\Omega)\,. 
\]
Then we set
\[
C^{k,\beta}_{b}({\mathrm{cl}}\Omega)\equiv
\{
u\in C^{k,\beta}({\mathrm{cl}}\Omega):\,
D^{\gamma}u\ {\mathrm{is\ bounded}}\ \forall\gamma\in {\mathbb{N}}^{n}\
{\mathrm{such\ that}}\ |\gamma|\leq k
\}\,,
\]
and we endow $C^{k,\beta}_{b}({\mathrm{cl}}\Omega)$ with its usual  norm
\[
\|u\|_{ C^{k,\beta}_{b}({\mathrm{cl}}\Omega) }\equiv
\sum_{|\gamma|\leq k}\sup_{x\in {\mathrm{cl}}\Omega }|D^{\gamma}u(x)|
+\sum_{|\gamma| = k}|D^{\gamma}u: {\mathrm{cl}}\Omega |_{\beta}
\qquad\forall u\in C^{k,\beta}_{b}({\mathrm{cl}}\Omega)\,,
\]
where $|D^{\gamma}u: {\mathrm{cl}}\Omega |_{\beta}$ denotes the $\beta$-H\"{o}lder constant of $D^{\gamma}u$.

Next we turn to periodic domains. If $\mathbb{I}$ is an arbitrary subset of ${\mathbb{R}}^{n}$  such that
 $\mathrm{cl}\mathbb{I}\subseteq Q$, then we set
 \begin{align}
&{\mathbb{S}} [\mathbb{I}]\equiv 
\bigcup_{z\in{\mathbb{Z}}^{n} }(qz+\mathbb{I})=q{\mathbb{Z}}^{n}+\mathbb{I}\,,\nonumber
\\
&{\mathbb{S}} [\mathbb{I}]^{-}\equiv {\mathbb{R}}^{n}\setminus{\mathrm{cl}}{\mathbb{S}} [\mathbb{I}]\,.\nonumber
\end{align}
We note that if $\mathbb{R}^n\setminus \mathrm{cl}\mathbb{I}$ is connected, then $\mathbb{S}[\mathbb{I}]^{-}$ is connected. 

Let $\mathbb{D}\subseteq \mathbb{R}^n$ be such that $qz+\mathbb{D}\subseteq \mathbb{D}$ for all $z \in \mathbb{Z}^n$. We say that a function $u$ from $\mathbb{D}$ to $\mathbb{R}$ is $q$--periodic if $u(x+qe_j)=u(x)$ for all $x \in \mathbb{D}$ and for all $j\in \{1,\dots,n\}$.

If $\mathbb{I}$ is an open subset of ${\mathbb{R}}^{n}$  such that ${\mathrm{cl}}\mathbb{I}\subseteq Q$ and if 
$k\in {\mathbb{N}}$, $\beta\in]0,1]$, then we set\[
C^{k}_{q}({\mathrm{cl}}{\mathbb{S}}[\mathbb{I}] )
\equiv\left\{
u\in C^{k}({\mathrm{cl}}{\mathbb{S}}[\mathbb{I}] ):\,
u\ {\mathrm{is}}\ q-{\mathrm{periodic}}
\right\}\,,
\]
which we regard as a Banach subspace of $C^{k}_{b}({\mathrm{cl}}{\mathbb{S}}[\mathbb{I}] )$, and 
\[
C^{k,\beta}_{q}({\mathrm{cl}}{\mathbb{S}}[\mathbb{I}] )
\equiv\left\{
u\in C^{k,\beta}({\mathrm{cl}}{\mathbb{S}}[\mathbb{I}] ):\,
u\ {\mathrm{is}}\ q-{\mathrm{periodic}}
\right\}\,,
\]
which we regard as a Banach subspace of $C^{k,\beta}_{b}({\mathrm{cl}}{\mathbb{S}}[\mathbb{I}] )$, and 
\[
C^{k}_{q}({\mathrm{cl}}{\mathbb{S}}[\mathbb{I}]^{-})
\equiv\left\{
u\in C^{k}({\mathrm{cl}}{\mathbb{S}}[\mathbb{I}]^{-}):\,
u\ {\mathrm{is}}\ q-{\mathrm{periodic}}
\right\}\,,
\]
which we regard as a Banach subspace of $C^{k}_{b}({\mathrm{cl}}{\mathbb{S}}[\mathbb{I}]^{-})$, and
\[
C^{k,\beta}_{q}({\mathrm{cl}}{\mathbb{S}}[\mathbb{I}]^{-} )
\equiv\left\{
u\in C^{k,\beta}({\mathrm{cl}}{\mathbb{S}}[\mathbb{I}]^{-} ):\,
u\ {\mathrm{is}}\ q-{\mathrm{periodic}}
\right\}\,,
\]
which we regard as a Banach subspace of $C^{k,\beta}_{b}({\mathrm{cl}}{\mathbb{S}}[\mathbb{I}]^{-})$.  
We denote by ${\mathcal{S}}( {\mathbb{R}}^{n})$ the Schwartz space of complex-valued rapidly decreasing functions. 

In the following Theorem, we introduce a periodic analog of the fundamental solution of the Laplace operator (cf.~\textit{e.g.}, Hasimoto \cite{Ha59}, Shcherbina \cite{Sh86}, Poulton, Botten, McPhedran, and Movchan \cite{PoBoMcMo99}, Ammari, Kang, and Touibi \cite{AmKaTo05}, Ammari and Kang \cite[p. 53]{AmKa07}, and \cite{LaMu10a}.)

\begin{thm}
\label{psper}
The generalized series
\begin{equation}
\label{psper2}
S_{n}^{q}
\equiv
\sum_{z\in {\mathbb{Z}}^{n}\setminus \{0\}} \frac{1}{-4\pi^2 |q^{-1}z|^2|Q|}E_{ 2\pi iq^{-1}z }
\end{equation}
defines a tempered distribution in ${\mathbb{R}}^{n}$ such that $S_{n}^{q}$ is $q$-periodic, \textit{i.e.}, 
\[
\tau_{q_{jj}e_{j}}S_{n}^{q}=S_{n}^{q}\qquad\forall j\in\{1,\dots,n\}\,,
\]
and such that 
\[
\Delta S_{n}^{q}=\sum_{z\in {\mathbb{Z}}^{n}}
\delta_{qz}-\frac{1}{|Q|}\,,
\]
where $\delta_{qz}$ denotes the Dirac measure with mass at $qz$, for all $z\in {\mathbb{Z}}^{n}$. Moreover, the following statements hold.
\begin{enumerate}
\item[(i)] $S_{n}^{q}$ is real analytic in ${\mathbb{R}}^{n}\setminus q{\mathbb{Z}}^{n}$.
\item[(ii)] $R_{n}^{q}\equiv S_{n}^{q}-S_{n}$ is real analytic in 
$({\mathbb{R}}^{n}\setminus q {\mathbb{Z}}^{n})\cup \{0\}$, and we have
\[
\Delta R_{n}^{q}
=
\sum_{z\in {\mathbb{Z}}^{n}\setminus\{0\} }\delta_{qz}
-
 \frac{1}{|Q|}
\,,
\]
in the sense of distributions.
\item[(iii)] $S_{n}^{q}\in L^{1}_{ {\mathrm{loc}} }({\mathbb{R}}^{n})$.
\item[(iv)] $S_{n}^{q}(x)=S_{n}^{q}(-x)$ for all $x \in \mathbb{R}^n \setminus q\mathbb{Z}^n$. 
\end{enumerate}
\end{thm}
{\bf Proof.} For the proof of (i), (ii), we refer for example to \cite{LaMu10a}, where an analog of a periodic fundamental solution for a second order strongly elliptic differential operator with constant coefficients has been constructed. We now consider statement (iii). As is well known, $S_n^q$ is a locally integrable complex-valued function (cf. \cite[\S 3]{LaMu10a}.) By the definition of $S_n^q$, and by the equality
\[
\overline{<E_{2\pi i q^{-1}z},\overline{\phi}>}=<E_{2\pi i q^{-1}(-z)},\phi> \qquad \forall \phi \in \mathcal{S}(\mathbb{R}^n)\, , \qquad \forall z \in \mathbb{Z}^n \setminus \{0\}\, ,
\]
and by the obvious identity
\[
\frac{1}{-4\pi^2 |-q^{-1}z|^2|Q|}=\frac{1}{-4\pi^2 |q^{-1}z|^2|Q|} \qquad \forall z \in \mathbb{Z}^n \setminus \{0\}\,,
\]
we can conclude that $S_{n}^{q}$ is actually a real-valued function. We now turn to the proof of (iv). By a straightforward verification based on \eqref{psper2}, we have
\[
\int_{\mathbb{R}^n}S_n^q(x)\phi(-x)\, dx=\int_{\mathbb{R}^n}S_n^q(x)\phi(x)\, dx \qquad \forall \phi \in \mathcal{S}(\mathbb{R}^n)\,,
\]
and thus $S_{n}^{q}(x)=S_{n}^{q}(-x)$ for all $x \in \mathbb{R}^n \setminus q\mathbb{Z}^n$.  Hence, the proof is complete \hfill $\Box$
\vspace{\baselineskip}

We now introduce the periodic double layer potential. Let  $\alpha\in]0,1[$, $m\in {\mathbb{N}}\setminus\{0\}$. Let $\mathbb{I}$ be a bounded connected open subset of ${\mathbb{R}}^{n}$ of class $C^{m,\alpha}$ such that 
${\mathbb{R}}^{n}\setminus{\mathrm{cl}}\mathbb{I}$ is connected and that ${\mathrm{cl}}\mathbb{I}\subseteq Q$. Let $S_{n}^{q}$ be as in Theorem  \ref{psper}. If $\mu\in C^{0,\alpha}(\partial\mathbb{I})$, we set
\[
w_{q}[\partial\mathbb{I},\mu](x)\equiv
-\int_{\partial\mathbb{I}}(DS_{n}^{q}(x-y) )\nu_{\mathbb{I}}(y)\mu(y)\,d\sigma_{y}
\qquad\forall x\in {\mathbb{R}}^{n}\,.
\]
In the following Theorem, we collect some properties of the periodic double layer potential.
\begin{thm}
\label{dperpot}
Let  $\alpha\in]0,1[$, $m\in {\mathbb{N}}\setminus\{0\}$. Let $\mathbb{I}$ be a bounded connected open subset of ${\mathbb{R}}^{n}$ of class $C^{m,\alpha}$ such that 
${\mathbb{R}}^{n}\setminus{\mathrm{cl}}\mathbb{I}$ is connected and that ${\mathrm{cl}}\mathbb{I}\subseteq Q$. Let $S_{n}^{q}$ be as in Theorem  \ref{psper}. Then the following statements hold.
\begin{enumerate}
\item[(i)] Let $\mu\in C^{0,\alpha}(\partial\mathbb{I})$. Then $w_{q}[\partial\mathbb{I},\mu]$ is $q$-periodic and
\[
\Delta(w_{q}[\partial\mathbb{I},\mu])(x)=0\qquad\forall x\in {\mathbb{R}}^{n}\setminus\partial{\mathbb{S}}[\mathbb{I}]\,.
\]
\item[(ii)] If $\mu\in C^{m,\alpha}(\partial\mathbb{I})$, then 
the restriction $w_{q}[\partial\mathbb{I},\mu]_{|{\mathbb{S}}[\mathbb{I}]}$ can be extended uniquely to an element $w^{+}_{ q }[\partial\mathbb{I},\mu]$ of $C_{q}^{m,\alpha}({\mathrm{cl}}{\mathbb{S}}[\mathbb{I}])$, and the  restriction $w_{q}[\partial\mathbb{I},\mu]_{|{\mathbb{S}}[\mathbb{I}]^{-}}$ can be extended uniquely to an element $w_{ q }^{-}[\partial\mathbb{I},\mu]$ of $C^{m,\alpha}_{q}({\mathrm{cl}}{\mathbb{S}}[\mathbb{I}]^{-})$, and we have
\begin{align}
&w^{\pm }_{q}[\partial\mathbb{I},\mu]=\pm\frac{1}{2}\mu+
w_{q}[\partial\mathbb{I},\mu]
\qquad{\mathrm{on}}\ \partial \mathbb{I} \,,\label{dperpot2a}
\\  
&(Dw^{+}_{q}[\partial\mathbb{I},\mu])\nu_{\mathbb{I}}-
(Dw^{-}_{q}[\partial\mathbb{I},\mu])\nu_{\mathbb{I}}=0 \ \  {\mathrm{on}}\ \partial\mathbb{I}\,.\label{dperpot2b}
\end{align}
\item[(iii)] The operator from  $C^{m,\alpha}(\partial\mathbb{I})$ to 
$C^{m,\alpha}_{q}({\mathrm{cl}}{\mathbb{S}}[\mathbb{I}])$ which takes $\mu$ to the function $
w_{q}^{+}[\partial\mathbb{I},\mu]$ is  continuous. The operator from   $C^{m,\alpha}(\partial\mathbb{I})$ to $C^{m,\alpha}_{q}({\mathrm{cl}}{\mathbb{S}}[\mathbb{I}]^{-})$ which takes $\mu$ to the function $w_{q}^{-}[\partial\mathbb{I},\mu]$ is continuous. 
\item[(iv)] The following equalities hold
\begin{align}
& w_{q}[\partial \mathbb{I},1](x)=\frac{1}{2}-\frac{|\mathbb{I}|}{|Q|} \qquad \forall x \in \partial \mathbb{S}[\mathbb{I}] \, ,\label{dperpot3a}\\
& w_{q}[\partial \mathbb{I},1](x)=1-\frac{|\mathbb{I}|}{|Q|} \qquad \forall x \in \mathbb{S}[\mathbb{I}] \, ,\label{dperpot3b}\\
& w_{q}[\partial \mathbb{I},1](x)=-\frac{|\mathbb{I}|}{|Q|} \qquad \forall x \in \mathbb{S}[\mathbb{I}]^{-} \, ,\label{dperpot3c}
\end{align}
where $|\mathbb{I}|$, $|Q|$ denote the $n$-dimensional measure of $\mathbb{I}$ and of $Q$, respectively.
\end{enumerate}
\end{thm}
{\bf Proof.} For the proof of statements (i), (ii), (iii), we refer for example to \cite{LaMu10a}. We now consider statement (iv). It clearly suffices to prove \eqref{dperpot3c}. Indeed, equalities \eqref{dperpot3a}, \eqref{dperpot3b} can be proved by exploiting \eqref{dperpot3c} and the jump relations of \eqref{dperpot2a}. By the periodicity of $w_q[\partial \mathbb{I},1]$, we can assume $x \in \mathrm{cl}Q \setminus \mathrm{cl}\mathbb{I}$. By the Green formula and Theorem \ref{psper}, we have
\[
-\int_{\partial\mathbb{I}}(DS_{n}^{q}(x-y) )\nu_{\mathbb{I}}(y)\,d\sigma_{y}=\int_{\mathbb{I}}\Delta_y(S_{n}^{q}(x-y) )\,dy=-\frac{|\mathbb{I}|}{|Q|}\,,
\]
and accordingly \eqref{dperpot3c} holds. Thus the proof is complete.
\hfill $\Box$
\vspace{\baselineskip}

Let  $\alpha\in]0,1[$, $m\in {\mathbb{N}}\setminus\{0\}$. If $\Omega$ is a bounded connected open subset of ${\mathbb{R}}^{n}$ of class $C^{m,\alpha}$, we find convenient to set
\[
C^{m,\alpha}(\partial \Omega)_0\equiv \left\{f \in C^{m,\alpha}(\partial \Omega)\colon \int_{\partial \Omega}f \, d\sigma=0\right\}\,.
\]
Then we have the following Proposition.
\begin{prop}\label{prop:linhom}
Let  $\alpha\in]0,1[$, $m\in {\mathbb{N}}\setminus\{0\}$. Let $\mathbb{I}$ be a bounded connected open subset of ${\mathbb{R}}^{n}$ of class $C^{m,\alpha}$ such that 
${\mathbb{R}}^{n}\setminus{\mathrm{cl}}\mathbb{I}$ is connected and that ${\mathrm{cl}}\mathbb{I}\subseteq Q$. Let $S_{n}^{q}$ be as in Theorem  \ref{psper}. Let $M[\cdot,\cdot]$ be the map from $C^{m,\alpha}(\partial \mathbb{I})_0\times \mathbb{R}$ to $C^{m,\alpha}(\partial \mathbb{I})$, defined by
\[
M[\mu,\xi](x)\equiv -\frac{1}{2}\mu(x)+w_{q}[\partial \mathbb{I},\mu](x)+\xi \qquad \forall x \in \partial \mathbb{I}\, ,
\]
for all $(\mu,\xi) \in C^{m,\alpha}(\partial \mathbb{I})_0\times \mathbb{R}$. Then $M[\cdot,\cdot]$ is a linear homeomorphism from $C^{m,\alpha}(\partial \mathbb{I})_0\times \mathbb{R}$ onto $C^{m,\alpha}(\partial \mathbb{I})$.
\end{prop}
 {\bf Proof.} By Theorem \ref{dperpot}, $M$ is continuous. As a consequence, by the Open Mapping Theorem, it suffices to prove that $M$ is a bijection. We first show that $M$ is injective. So let $(\mu,\xi) \in C^{m,\alpha}(\partial \mathbb{I})_0\times \mathbb{R}$ be such that $M[\mu,\xi]=0$. Then,
 \[
-\frac{1}{2}\mu(x)+w_{q}[\partial \mathbb{I},\mu](x)=-\xi \qquad \forall x \in \partial \mathbb{I}\, , 
\]
 and thus, by Proposition \ref{prop:comphom} of the Appendix, $\mu$ must be constant. Since $\int_{\partial \mathbb{I}} \mu\, d\sigma=0$, then $\mu=0$, and so also $\xi=0$. Hence $M$ is injective. It remains to prove that $M$ is surjective. So let $g \in C^{m,\alpha}(\partial \mathbb{I})$. We need to prove that there exists a pair $(\mu,\xi) \in C^{m,\alpha}(\partial \mathbb{I})_0\times \mathbb{R}$ such that $M[\mu,\xi]=g$. By Proposition \ref{prop:comphom} of the Appendix, there exists a (unique) $\tilde{\mu} \in C^{m,\alpha}(\partial \mathbb{I})$ such that
 \[
-\frac{1}{2}\tilde{\mu}(x)+w_{q}[\partial \mathbb{I},\tilde{\mu}](x)=g(x) \qquad \forall x \in \partial \mathbb{I}\, .
 \]
 Accordingly, if we set
 \begin{align}
 &\mu(x)\equiv \tilde{\mu}(x)-\frac{1}{\int_{\partial \mathbb{I}}d\sigma}\int_{\partial \mathbb{I}}\tilde{\mu}\, d\sigma \qquad \forall x \in \partial \mathbb{I}\, ,\nonumber \\
 &\xi \equiv -\frac{|\mathbb{I}|}{|Q|}\frac{1}{\int_{\partial \mathbb{I}}d\sigma}\int_{\partial \mathbb{I}}\tilde{\mu}\, d\sigma \, ,\nonumber
 \end{align}
 where $|\mathbb{I}|$, $|Q|$ denote the $n$-dimensional measure of $\mathbb{I}$ and of $Q$, respectively, then clearly $(\mu,\xi) \in C^{m,\alpha}(\partial \mathbb{I})_0\times \mathbb{R}$ and $M[\mu,\xi]=g$. Therefore, $M$ is bijective, and the proof is complete.
 \hfill $\Box$
\vspace{\baselineskip}

In the following Proposition, we show that a periodic Dirichlet boundary value problem in the perforated domain $\mathbb{S}[\mathbb{I}]^{-}$ has a unique solution in $C^{m,\alpha}_q(\mathrm{cl}\mathbb{S}[\mathbb{I}]^{-})$, which can be represented as the sum of a periodic double layer potential and a costant.
\begin{prop}\label{prop:Dirsol}
Let  $\alpha\in]0,1[$, $m\in {\mathbb{N}}\setminus\{0\}$. Let $\mathbb{I}$ be a bounded connected open subset of ${\mathbb{R}}^{n}$ of class $C^{m,\alpha}$ such that 
${\mathbb{R}}^{n}\setminus{\mathrm{cl}}\mathbb{I}$ is connected and that ${\mathrm{cl}}\mathbb{I}\subseteq Q$. Let $S_{n}^{q}$ be as in Theorem  \ref{psper}. Let $\Gamma \in C^{m,\alpha}(\partial \mathbb{I})$. Then the following boundary value problem
 \begin{equation}\label{bvp:Dir}
 \left \lbrace 
 \begin{array}{ll}
 \Delta u (x)= 0 & \textrm{$\forall x \in {\mathbb{S}} [\mathbb{I}]^{-}$}\,, \\
u(x+qe_i) =u(x) &  \textrm{$\forall x \in \mathrm{cl} {\mathbb{S}} [\mathbb{I}]^{-}\,, \quad \forall i \in \{1,\dots,n\}$}\,, \\
u(x)=\Gamma(x)& \textrm{$\forall x \in \partial \mathbb{I}$}\,,
 \end{array}
 \right.
 \end{equation}
 has a unique solution $u \in C^{m,\alpha}_q(\mathrm{cl} \mathbb{S}[\mathbb{I}]^{-})$. Moreover,
\begin{equation}\label{eq:Dirsol1}
u(x)=w_{q}^{-}[\partial \mathbb{I},\mu](x)+\xi \qquad \forall x \in \mathrm{cl}\mathbb{S}[\mathbb{I}]^{-}\,,
\end{equation}
where $(\mu,\xi)$ is the unique solution in $C^{m,\alpha}(\partial \mathbb{I})_0\times \mathbb{R}$ of the following integral equation
\begin{equation}\label{eq:Dirsol2}
\Gamma(x)=-\frac{1}{2}\mu(x)+w_{q}[\partial \mathbb{I},\mu](x)+\xi \qquad \forall x \in \partial \mathbb{I}\, .
\end{equation}
\end{prop}
 {\bf Proof.} We first note that Proposition \ref{prop:maxprin} of the Appendix implies that problem \eqref{bvp:Dir} has at most one solution. As a consequence, we need to prove that the function defined by \eqref{eq:Dirsol1} solves problem \eqref{bvp:Dir}. By Proposition \ref{prop:linhom}, there exists a unique solution $(\mu,\xi)\in C^{m,\alpha}(\partial \mathbb{I})_0\times \mathbb{R}$ of \eqref{eq:Dirsol2}. Then, by Theorem \ref{dperpot} and equation \eqref{eq:Dirsol2}, the function defined by \eqref{eq:Dirsol1} is a periodic harmonic function satisfying the third condition of \eqref{bvp:Dir}, and thus a solution of problem \eqref{bvp:Dir}. \hfill $\Box$
\vspace{\baselineskip}

\begin{rem}
Let the assumptions of Proposition \ref{prop:Dirsol} hold. We note that we proved, in particular, that the solution of boundary value problem \eqref{bvp:Dir} can be represented as the sum of a periodic double layer potential and a constant. However, we observe that we could also represent the solution of problem \eqref{bvp:Dir} as a periodic double layer potential (cf. Proposition \ref{prop:comphom} of the Appendix.) On the other hand, for the analysis of \eqref{bvp:Direps} around the degenerate value $(\epsilon,g)=(0,g_0)$, it will be preferable to exploit the representation formula of Proposition \ref{prop:Dirsol}.
\end{rem}

\section{Formulation of the problem in terms of integral equations}\label{form}

We now provide a formulation of problem \eqref{bvp:Direps} in terms of an integral equation. We shall consider the following assumptions for some $\alpha \in ]0,1[$ and for some natural $m \geq 1$.
\begin{equation}\label{ass} 
\begin{split}
&\text{Let $\Omega$ be a bounded connected open subset of ${\mathbb{R}}^{n}$ of class $C^{m,\alpha}$ such that 
${\mathbb{R}}^{n}\setminus{\mathrm{cl}}\Omega$ is connected and that $0 \in\Omega$}.\\
&
\text{Let $p \in Q$}. 
\end{split}
\end{equation}
If $\epsilon \in \mathbb{R}$, we set
\[
\Omega_\epsilon \equiv p+\epsilon\Omega \,.
\]
Now let 
\begin{equation}\label{eps0}
\epsilon_0\equiv \sup \Bigl \{ \theta \in ]0,+\infty[\colon \mathrm{cl}\Omega_\epsilon \subseteq Q\, , \forall \epsilon \in ]-\theta,\theta[ \Bigr\}\,.
\end{equation}
A simple topological argument shows that if \eqref{ass} holds, then $\mathbb{S}[\Omega_\epsilon]^{-}$ is connected, for all $\epsilon \in ]-\epsilon_0,\epsilon_0[$. We also note that
\[
\nu_{\Omega_\epsilon}(p+\epsilon t)=\mathrm{sgn}(\epsilon)\nu_{\Omega}(t) \qquad \forall t \in \partial \Omega\, ,
\]
for all $\epsilon \in ]-\epsilon_0,\epsilon_0[\setminus \{0\}$, where $\mathrm{sgn}(\epsilon)=1$ if $\epsilon>0$, $\mathrm{sgn}(\epsilon)=-1$ if $\epsilon<0$. Then we shall consider the following assumption.
\begin{equation}\label{assg}
\text{Let $g_0 \in C^{m,\alpha}(\partial \Omega)$.}
\end{equation}
If $(\epsilon,g) \in ]0,\epsilon_0[\times C^{m,\alpha}(\partial \Omega)$, we shall convert our boundary value problem \eqref{bvp:Direps} into an integral equation. We could exploit Proposition \ref{prop:Dirsol}, with $\mathbb{I}$ replaced by $\Omega_\epsilon$, but we note that the integral equation and the corresponding integral representation of the solution include integrations on the $\epsilon$-dependent domain $\partial \Omega_\epsilon$. In order to get rid of such a dependence, we shall introduce the following Lemma, in which we properly rescale the unknown density.

\begin{lem}\label{lem:equiv}
Let $\alpha \in ]0,1[$. Let $m \in \mathbb{N}\setminus \{0\}$. Let \eqref{ass}-\eqref{assg} hold. Let $S_{n}^{q}$, $R_{n}^{q}$ be as in Theorem \ref{psper}. Let $(\epsilon,g) \in ]0,\epsilon_0[\times C^{m,\alpha}(\partial \Omega)$. Then a pair $(\theta,\xi) \in C^{m,\alpha}(\partial \Omega)_0\times \mathbb{R}$ solves equation
\begin{equation}\label{eq:equiv1}
-\frac{1}{2}\theta(t)-\int_{\partial\Omega}(DS_{n}(t-s) )\nu_{\Omega}(s)\theta(s)\,d\sigma_{s}-\epsilon^{n-1}\int_{\partial\Omega}(DR_{n}^{q}(\epsilon(t-s)) )\nu_{\Omega}(s)\theta(s)\,d\sigma_{s}+\xi=g(t)\qquad \forall t \in \partial \Omega\, ,
\end{equation}
if and only if the pair $(\mu,\xi) \in C^{m,\alpha}(\partial \Omega_\epsilon)_0\times \mathbb{R}$, with $\mu$ delivered by
\begin{equation}\label{eq:equiv3}
\mu(x)\equiv \theta\Bigl(\frac{1}{\epsilon}(x-p)\Bigr) \qquad \forall x \in \partial \Omega_\epsilon\, ,
\end{equation}
 solves equation 
\begin{equation}\label{eq:equiv2}
-\frac{1}{2}\mu(x)+w_{q}[\partial \Omega_\epsilon,\mu](x)+\xi =\Gamma(x)\qquad \forall x \in \partial \Omega_\epsilon\, ,
\end{equation}
where
\[
\Gamma(x)\equiv g\Bigl(\frac{1}{\epsilon}(x-p)\Bigr) \qquad \forall x \in \partial \Omega_\epsilon\, .
\]
Moreover, equation \eqref{eq:equiv1} has a unique solution in $C^{m,\alpha}(\partial \Omega)_0\times \mathbb{R}$.
\end{lem}
 {\bf Proof.} The equivalence of equation \eqref{eq:equiv1} in the unknown $(\theta,\xi)$ and equation \eqref{eq:equiv2} in the unknown $(\mu,\xi)$, with $\mu$ delivered by \eqref{eq:equiv3}, is a straightforward consequence of the Theorem of change of variables in integrals. Then the existence and uniqueness of a solution in $C^{m,\alpha}(\partial \Omega)_0\times \mathbb{R}$ of equation \eqref{eq:equiv1}, follows from Proposition \ref{prop:linhom} applied to equation \eqref{eq:equiv2}, and from the equivalence of equations \eqref{eq:equiv1}, \eqref{eq:equiv2}. \hfill $\Box$
\vspace{\baselineskip}

In the following Lemma, we study equation \eqref{eq:equiv1}, when $(\epsilon,g)=(0,g_0)$.

\begin{lem}\label{lem:limeq}
Let $\alpha \in ]0,1[$. Let $m \in \mathbb{N}\setminus \{0\}$. Let \eqref{ass}-\eqref{assg} hold. Let $\tau_0$ be the unique solution in $C^{m-1,\alpha}(\partial \Omega)$ of the following problem
\begin{equation}\label{eq:limeq0}
  \left \lbrace 
 \begin{array}{ll}
-\frac{1}{2}\tau(t)+\int_{\partial\Omega}(DS_{n}(t-s) )\nu_{\Omega}(t)\tau(s)\,d\sigma_{s}=0 \qquad \forall t \in \partial \Omega\, ,\\
\int_{\partial \Omega}\tau \, d\sigma=1\, . & 
 \end{array}
 \right.
\end{equation}
Then equation
\begin{equation}\label{eq:limeq1}
-\frac{1}{2}\theta(t)-\int_{\partial\Omega}(DS_{n}(t-s) )\nu_{\Omega}(s)\theta(s)\,d\sigma_{s}+\xi=g_0(t)\qquad \forall t \in \partial \Omega\, ,
\end{equation}
which we call the \emph{limiting equation}, has a unique solution in $C^{m,\alpha}(\partial \Omega)_0\times \mathbb{R}$, which we denote by $(\tilde{\theta},\tilde{\xi})$. Moreover, 
\begin{equation}\label{eq:limeq1a}
\tilde{\xi}=\int_{\partial \Omega}g_0\tau_0\, d\sigma\,,
\end{equation}
and the function $\tilde{u} \in C^{m,\alpha}(\mathbb{R}^n \setminus \Omega)$, defined by
\begin{equation}\label{eq:limeq2}
\tilde{u}(t)\equiv -\int_{\partial \Omega}(DS_n(t-s))\nu_{\Omega}(s)\tilde{\theta}(s)d\sigma_s \qquad \forall t \in \mathbb{R}^n \setminus \mathrm{cl}\Omega\, ,
\end{equation}
and extended by continuity to $\mathbb{R}^{n}\setminus \Omega$, is the unique solution in $C^{m,\alpha}(\mathbb{R}^n\setminus \Omega)$ of the following problem
\begin{equation}\label{eq:limeq3}
 \left \lbrace 
 \begin{array}{ll}
 \Delta u (t)= 0 & \textrm{$\forall t \in \mathbb{R}^n \setminus \mathrm{cl}\Omega$}\,, \\
u(t)=g_0(t)-\int_{\partial \Omega}g_0\tau_0\, d\sigma&  \textrm{$\forall t \in \partial \Omega$}\,, \\
\lim_{t\to \infty}u(t)=0\,.&
 \end{array}
 \right.
\end{equation}
\end{lem}
 {\bf Proof.} We first note that the unique solvability of problem \eqref{eq:limeq0} in the class of continuous functions follows by classical potential theory (cf.~\textit{e.g.}, Folland \cite[Ch.~3]{Fo95}.) For the $C^{m-1,\alpha}$ regularity of the solution, we refer, \textit{e.g.}, to Lanza \cite[Appendix A]{La07a}. By Proposition \ref{prop:classpot} of the Appendix, equation \eqref{eq:limeq1} has a unique solution in $C^{m,\alpha}(\partial \Omega)_0\times \mathbb{R}$. Moreover, as is well known, if $\psi \in C^{m,\alpha}(\partial \Omega)$, then
\[
\psi \in\Bigl\{-\frac{1}{2}\theta(\cdot)-\int_{\partial\Omega}(DS_{n}(\cdot-s) )\nu_{\Omega}(s)\theta(s)\,d\sigma_{s} \colon \theta \in C^{m,\alpha}(\partial \Omega)\Bigr\}
\] 
if and only if
\[
\int_{\partial \Omega}\psi \tau_0\, d\sigma=0\,,
\]
and thus $\tilde{\xi}$ must be delivered by equality \eqref{eq:limeq1a} (cf.~\textit{e.g.}, Folland \cite[Ch.~3]{Fo95} and Lanza \cite[Appendix A]{La07a}.) Then by classical potential theory, the function defined by \eqref{eq:limeq2} and extended by continuity to $\mathbb{R}^n \setminus \Omega$ solves problem \eqref{eq:limeq3}, which has at most one solution (cf.~\textit{e.g.}, Folland \cite[Ch.~3]{Fo95}, Miranda \cite{Mi65}, Dalla Riva and Lanza \cite[Theorem 3.1]{DaLa10a}, Lanza and Rossi \cite[Theorem 3.1]{LaRo04}.)\hfill $\Box$
\vspace{\baselineskip}

We are now ready to analyse equation \eqref{eq:equiv1} around the degenerate case $(\epsilon,g)=(0,g_0)$.
We find convenient to introduce the following abbreviation. We set
\[
\mathcal{X}_{m,\alpha}\equiv C^{m,\alpha}(\partial \Omega)_0\times \mathbb{R}\,.
\]
Then we have the following.
\begin{prop}\label{prop:Lambda}
Let $\alpha \in ]0,1[$. Let $m \in \mathbb{N}\setminus \{0\}$. Let \eqref{ass}-\eqref{assg} hold. Let $S_{n}^{q}$, $R_{n}^{q}$ be as in Theorem \ref{psper}. Let $\Lambda$ be the map from $]-\epsilon_0,\epsilon_0[\times C^{m,\alpha}(\partial \Omega)\times \mathcal{X}_{m,\alpha}$ to $C^{m,\alpha}(\partial \Omega)$, defined by
\[
\Lambda[\epsilon,g,\theta,\xi](t)\equiv -\frac{1}{2}\theta(t)-\int_{\partial\Omega}(DS_{n}(t-s) )\nu_{\Omega}(s)\theta(s)\,d\sigma_{s}-\epsilon^{n-1}\int_{\partial\Omega}(DR_{n}^{q}(\epsilon(t-s)) )\nu_{\Omega}(s)\theta(s)\,d\sigma_{s} +\xi-g(t)\qquad \forall t \in \partial \Omega\, ,
\]
for all $(\epsilon,g,\theta,\xi) \in ]-\epsilon_0,\epsilon_0[\times C^{m,\alpha}(\partial \Omega)\times \mathcal{X}_{m,\alpha}$. Then the following statements hold.
\begin{enumerate}
\item[(i)] Let $(\epsilon,g) \in ]0,\epsilon_0[\times C^{m,\alpha}(\partial \Omega)$. Then equation
\[
\Lambda[\epsilon,g,\theta,\xi]=0
\]
has a unique solution in $\mathcal{X}_{m,\alpha}$, which we denote by $(\hat{\theta}[\epsilon,g],\hat{\xi}[\epsilon,g])$ (cf. Lemma \ref{lem:equiv}.)
\item[(ii)] Equation
\[
\Lambda[0,g_0,\theta,\xi]=0
\]
has a unique solution in $\mathcal{X}_{m,\alpha}$, which we denote by $(\hat{\theta}[0,g_0],\hat{\xi}[0,g_0])$. Moreover,  $(\hat{\theta}[0,g_0],\hat{\xi}[0,g_0])=(\tilde{\theta},\tilde{\xi})$ (cf. Lemma \ref{lem:limeq}.)
\item[(iii)] $\Lambda[\cdot,\cdot,\cdot,\cdot]$ is a real analytic map from $]-\epsilon_0,\epsilon_0[\times C^{m,\alpha}(\partial \Omega)\times \mathcal{X}_{m,\alpha}$ to $C^{m,\alpha}(\partial \Omega)$. Moreover, the differential $\partial_{(\theta,\xi)}\Lambda[0,g_0,\hat{\theta}[0,g_0],\hat{\xi}[0,g_0]]$ of $\Lambda$ at $(0,g_0,\hat{\theta}[0,g_0],\hat{\xi}[0,g_0])$ with respect to the variables $(\theta,\xi)$ is a linear homeomorphism from $\mathcal{X}_{m,\alpha}$ onto $C^{m,\alpha}(\partial \Omega)$.
\item[(iv)] There exist $\epsilon_1 \in ]0,\epsilon_0]$, an open neighbourhood $\mathcal{U}$ of $g_0$ in $C^{m,\alpha}(\partial \Omega)$, and a real analytic map $(\Theta[\cdot,\cdot],\Xi[\cdot,\cdot])$ from $]-\epsilon_1,\epsilon_1[\times \mathcal{U}$ to $\mathcal{X}_{m,\alpha}$, such that
\begin{align}
&(\Theta[\epsilon,g],\Xi[\epsilon,g])=(\hat{\theta}[\epsilon,g],\hat{\xi}[\epsilon,g]) \qquad \forall (\epsilon,g) \in ]0,\epsilon_1[\times \mathcal{U}\, ,\nonumber \\
&(\Theta[0,g_0], \Xi[0,g_0])=(\tilde{\theta},\tilde{\xi})\, .\nonumber
\end{align}
\end{enumerate}
\end{prop}
 {\bf Proof.} Statements (i), (ii) are immediate consequences of Lemmas \ref{lem:equiv}, \ref{lem:limeq}. We now consider statement (iii). We first introduce some notation. For each $j \in \{1,\dots,n\}$, we denote by $R_j[\cdot,\cdot]$ the map from $]-\epsilon_0,\epsilon_0[\times L^1(\partial \Omega)$ to $C^{m,\alpha}(\mathrm{cl} \Omega)$, defined by
\[
R_j[\epsilon,f](t)\equiv \int_{\partial\Omega}(D_{x_j}R_{n}^{q}(\epsilon(t-s)) )f(s)\,d\sigma_{s}\qquad \forall t \in \mathrm{cl} \Omega\, ,
\]
for all $(\epsilon,f) \in ]-\epsilon_0,\epsilon_0[\times L^1(\partial \Omega)$. By classical potential theory and standard calculus in Banach spaces, we note that the map from $C^{m,\alpha}(\partial \Omega)\times C^{m,\alpha}(\partial \Omega)_0\times \mathbb{R}$ to $C^{m,\alpha}(\partial \Omega)$, which takes $(g,\theta,\xi)$ to the function 
\[
-\frac{1}{2}\theta(t)-\int_{\partial\Omega}(DS_{n}(t-s) )\nu_{\Omega}(s)\theta(s)\,d\sigma_{s}+\xi-g(t) 
\]
of the variable $t \in \partial \Omega$, is linear and continuous, and thus real analytic (cf.~\textit{e.g.}, Miranda \cite{Mi65}, Dalla Riva and Lanza \cite[Theorem 3.1]{DaLa10a}, Lanza and Rossi \cite[Theorem 3.1]{LaRo04}.) Then, in order to prove the real analyticity of $\Lambda[\cdot,\cdot,\cdot,\cdot]$ in $]-\epsilon_0,\epsilon_0[\times C^{m,\alpha}(\partial \Omega)\times \mathcal{X}_{m,\alpha}$, it clearly suffices to show that $R_j[\cdot,\cdot]$ is real analytic in $]-\epsilon_0,\epsilon_0[\times L^1(\partial \Omega)$ for each $j \in \{1,\dots,n\}$. Indeed, if $R_j[\cdot,\cdot]$ is real analytic from $]-\epsilon_0,\epsilon_0[\times L^1(\partial \Omega)$ to $C^{m,\alpha}(\mathrm{cl} \Omega)$ for all $j \in \{1,\dots,n\}$, then by the continuity of the linear map from $C^{m,\alpha}(\partial \Omega)_0$ to $L^1(\partial \Omega)$ which takes $\theta$ to $(\nu_{\Omega})_j\theta$, and by the continuity of the trace operator from $C^{m,\alpha}(\mathrm{cl}\Omega)$ to $C^{m,\alpha}(\partial \Omega)$, we can deduce the analyticity of  $\Lambda[\cdot,\cdot,\cdot,\cdot]$ in $]-\epsilon_0,\epsilon_0[\times C^{m,\alpha}(\partial \Omega)\times \mathcal{X}_{m,\alpha}$. Now let $\mathrm{id}_{\partial \Omega}$ and $\mathrm{id}_{\mathrm{cl}\Omega}$ denote the identity on $\partial \Omega$ and on $\mathrm{cl}\Omega$, respectively. Then we note that the map from $]-\epsilon_0,\epsilon_0[$ to $C^{m,\alpha}(\mathrm{cl}\Omega,\mathbb{R}^n)$ which takes $\epsilon$ to $\epsilon \mathrm{id}_{\mathrm{cl}\Omega}$, and the map from $]-\epsilon_0,\epsilon_0[$ to $C^{m,\alpha}(\partial\Omega,\mathbb{R}^n)$ which takes $\epsilon$ to $\epsilon \mathrm{id}_{\partial\Omega}$ are real analytic. Moreover,
\[
\epsilon \mathrm{cl}\Omega-\epsilon \partial \Omega \subseteq (\mathbb{R}^n \setminus q\mathbb{Z}^n)\cup \{0\} \qquad \forall \epsilon \in ]-\epsilon_0,\epsilon_0[\,.
\]
Then by the real analyticity of $D_{x_j}R_n^q$ in $(\mathbb{R}^n\setminus q\mathbb{Z}^n) \cup \{0\}$ and by Proposition \ref{prop:inop} (ii) of the Appendix, $R_j[\cdot,\cdot]$ is real analytic in $]-\epsilon_0,\epsilon_0[\times L^1(\partial \Omega)$, for each $j \in \{1,\dots,n\}$. Hence, $\Lambda[\cdot,\cdot,\cdot,\cdot]$ is real analytic in $]-\epsilon_0,\epsilon_0[\times C^{m,\alpha}(\partial \Omega)\times \mathcal{X}_{m,\alpha}$. By standard calculus in Banach space, the differential $\partial_{(\theta,\xi)}\Lambda[0,g_0,\hat{\theta}[0,g_0],\hat{\xi}[0,g_0]]$ of $\Lambda$ at $(0,g_0,\hat{\theta}[0,g_0],\hat{\xi}[0,g_0])$ with respect to $(\theta,\xi)$ is delivered by the following formula:
\[
\partial_{(\theta,\xi)}\Lambda[0,g_0,\hat{\theta}[0,g_0],\hat{\xi}[0,g_0]](\psi,\rho)(t)=-\frac{1}{2}\psi(t)-\int_{\partial\Omega}(DS_{n}(t-s) )\nu_{\Omega}(s)\psi(s)\,d\sigma_{s}+\rho\qquad \forall t \in \partial \Omega\, ,
\]
for all $(\psi,\rho) \in \mathcal{X}_{m,\alpha}$. Accordingly, by Proposition \ref{prop:classpot} of the Appendix, $\partial_{(\theta,\xi)}\Lambda[0,g_0,\hat{\theta}[0,g_0],\hat{\xi}[0,g_0]]$ is a linear homeomorphism from $\mathcal{X}_{m,\alpha}$ onto $C^{m,\alpha}(\partial \Omega)$, and the proof of (iii) is complete. Finally, statement (iv) is an immediate consequence of statement (iii) and of the Implicit Function Theorem for real analytic maps in Banach spaces (cf.~\textit{e.g.}, Prodi and Ambrosetti \cite[Theorem 11.6]{PrAm73}, Deimling \cite[Theorem 15.3]{De85}.)\hfill $\Box$
\vspace{\baselineskip}

\begin{rem}\label{rem:sol}
Let the assumptions of Proposition \ref{prop:Lambda} hold. Let $\epsilon_1$, $\mathcal{U}$, $(\Theta[\cdot,\cdot],\Xi[\cdot,\cdot])$ be as in Proposition \ref{prop:Lambda} (iv). Then, by the rule of change of variables in integrals, by Propositions \ref{prop:Dirsol}, \ref{prop:Lambda}, and by Lemma \ref{lem:equiv}, we have
\[
u[\epsilon,g](x)=-\epsilon^{n-1}\int_{\partial \Omega}(DS_{n}^{q}(x-p-\epsilon s))\nu_{\Omega}(s)\Theta[\epsilon,g](s)\, d\sigma_s+\Xi[\epsilon,g] \qquad \forall x \in \mathbb{S}[\Omega_\epsilon]^{-}\, ,
\]
for all $(\epsilon,g) \in ]0,\epsilon_1[\times \mathcal{U}$.
\end{rem}

\section{A functional analytic representation Theorem for the solution and its energy integral}\label{rep}

The following statement shows that suitable restrictions of $u[\epsilon,g](\cdot)$ can be continued real analytically for negative values of $\epsilon$.
\begin{thm}\label{thm:rep}
Let $\alpha \in ]0,1[$. Let $m \in \mathbb{N}\setminus \{0\}$. Let \eqref{ass}-\eqref{assg} hold. Let $\tilde{u}$ be as in Lemma \ref{lem:limeq}. Let $\epsilon_1$, $\mathcal{U}$, $\Xi[\cdot,\cdot]$ be as in Proposition \ref{prop:Lambda} (iv). Then the following statements hold.
\begin{enumerate}
\item[(i)] Let $V$ be a bounded open subset of $\mathbb{R}^n$ such that $\mathrm{cl}V \subseteq \mathbb{R}^n \setminus (p+q\mathbb{Z}^n)$. Let $r \in \mathbb{N}$. Then there exist $\epsilon_2 \in ]0,\epsilon_1]$ and a real analytic map $U$ from $]-\epsilon_2,\epsilon_2[\times \mathcal{U}$ to $C^{r}(\mathrm{cl}V)$ such that the following statements hold.
\begin{enumerate}
\item[(j)] $\mathrm{cl}V \subseteq \mathbb{S}[\Omega_\epsilon]^{-}$ for all $\epsilon \in ]-\epsilon_2,\epsilon_2[$.
\item[(jj)]
\begin{equation}\label{eq:rep}
u[\epsilon,g](x)=\epsilon^{n-1}U[\epsilon,g](x)+\Xi[\epsilon,g] \qquad \forall x \in \mathrm{cl}V\,,
\end{equation}
for all $(\epsilon,g) \in ]0,\epsilon_2[\times \mathcal{U}$. Moreover,
\begin{equation}\label{eq:repbis}
U[0,g_0](x)=DS_{n}^{q}(x-p)\int_{\partial \Omega}\nu_{\Omega}(s)g_0(s)\, d\sigma_s - DS_{n}^{q}(x-p)\int_{\partial \Omega}s\frac{\partial \tilde{u}}{\partial \nu_{\Omega}}(s)\, d\sigma_s \qquad \forall x \in \mathrm{cl}V\, .
\end{equation}
\end{enumerate}
\item[(ii)] Let $\widetilde{V}$ be a bounded open subset of $\mathbb{R}^n \setminus \mathrm{cl}\Omega$. Then there exist $\tilde{\epsilon}_2 \in ]0,\epsilon_1]$ and a real analytic map $\widetilde{U}$ from $]-\tilde{\epsilon}_2,\tilde{\epsilon}_2[\times \mathcal{U}$ to $C^{m,\alpha}(\mathrm{cl}\widetilde{V})$ such that the following statements hold.
\begin{enumerate}
\item[(j')] $p+\epsilon\mathrm{cl}\widetilde{V}\subseteq Q\setminus \Omega_\epsilon$ for all $\epsilon \in ]-\tilde{\epsilon}_2,\tilde{\epsilon}_2[\setminus \{0\}$.
\item[(jj')]
\begin{equation}\label{eq:repa}
u[\epsilon,g](p+\epsilon t)=\widetilde{U}[\epsilon,g](t)+\Xi[\epsilon,g] \qquad \forall t \in \mathrm{cl}\widetilde{V}\,,
\end{equation}
for all $(\epsilon,g) \in ]0,\tilde{\epsilon}_2[\times \mathcal{U}$. Moreover,
\[
\widetilde{U}[0,g_0](t)=\tilde{u}(t) \qquad \forall t \in \mathrm{cl}\widetilde{V}\,.
\]
\end{enumerate}
\end{enumerate}
\end{thm}
{\bf Proof.}  Let $S_{n}^{q}$, $R_{n}^{q}$ be as in Theorem \ref{psper}. Let $(\Theta[\cdot,\cdot],\Xi[\cdot,\cdot])$ be as in Proposition \ref{prop:Lambda} (iv). We start by proving (i). By taking $\epsilon_2 \in ]0,\epsilon_1]$ small enough, we can clearly assume that (j) holds. Consider now (jj). If $(\epsilon,g) \in ]0,\epsilon_2[\times \mathcal{U}$, then by Remark \ref{rem:sol} we have
\[
u[\epsilon,g](x)=-\epsilon^{n-1}\int_{\partial \Omega}(DS_{n}^{q}(x-p-\epsilon s))\nu_{\Omega}(s)\Theta[\epsilon,g](s)\, d\sigma_s+\Xi[\epsilon,g] \qquad \forall x \in \mathrm{cl}V\, .
\]
Thus it is natural to set
\[
U[\epsilon,g](x)\equiv-\int_{\partial \Omega}(DS_{n}^{q}(x-p-\epsilon s))\nu_{\Omega}(s)\Theta[\epsilon,g](s)\, d\sigma_s \qquad \forall x \in \mathrm{cl}V\, ,
\]
for all $(\epsilon,g) \in ]-\epsilon_2,\epsilon_2[\times \mathcal{U}$. Then we note that
\[
\mathrm{cl}V-p-\epsilon \partial \Omega \subseteq \mathbb{R}^n \setminus q\mathbb{Z}^n \qquad \forall \epsilon \in ]-\epsilon_2,\epsilon_2[\,.
\]
As a consequence, by the real analyticity of $S_n^q$ in $\mathbb{R}^n \setminus q\mathbb{Z}^n$, and by the real analyticity of $\Theta[\cdot,\cdot]$ from $]-\epsilon_1,\epsilon_1[\times \mathcal{U}$ to $C^{m,\alpha}(\partial \Omega)_0$, and by Proposition \ref{prop:inop} (i) of the Appendix, we can conclude that $U$ is real analytic from $]-\epsilon_2,\epsilon_2[\times \mathcal{U}$ to $C^{r}(\mathrm{cl}V)$. By the definition of $U$, equality \eqref{eq:rep} holds for all $(\epsilon,g) \in ]0,\epsilon_2[\times \mathcal{U}$.
Next we turn to prove formula \eqref{eq:repbis}. First we note that
\[
U[0,g_0](x)=-DS_n^q(x-p)\int_{\partial \Omega}\nu_{\Omega}(s)\Theta[0,g_0](s)\, d\sigma_s \qquad \forall x \in \mathrm{cl}V\, .
\]
Proposition \ref{prop:Lambda} (iv) implies that $\Theta[0,g_0]=\tilde{\theta}$, where $\tilde{\theta}$ is as in Lemma \ref{lem:limeq}. Then we set
\[
w(t)\equiv -\int_{\partial \Omega}\bigl(DS_n(t-s)\bigr)\nu_{\Omega}(s)\tilde{\theta}(s)\, d\sigma_s \qquad \forall t \in \mathbb{R}^n\, .
\]
As is well known, $w_{|\Omega}$ admits a continuous extension to $\mathrm{cl}\Omega$, which we denote by $w^+$, and $w_{|\mathbb{R}^n \setminus \mathrm{cl}\Omega}$ admits a continuous extension to $\mathbb{R}^n \setminus \Omega$, which we denote by $w^-$. Moreover, $w^+ \in C^{m,\alpha}(\mathrm{cl}\Omega)$ and $w^- \in C^{m,\alpha}(\mathbb{R}^n \setminus \Omega)$ (cf. \textit{e.g.}, Lanza and Rossi \cite[Thm.~3.1]{LaRo04}.) Clearly, $w^- = \tilde{u}$.
Then we fix
 $j \in \{1,\dots,n\}$. By classical potential theory, we have
\[
\int_{\partial \Omega}\bigl(\nu_{\Omega}(s)\bigr)_j \tilde{\theta}(s)\, d\sigma_s=\int_{\partial \Omega}\bigl(\nu_{\Omega}(s)\bigr)_jw^+(s)\, d\sigma_s-\int_{\partial \Omega}\bigl(\nu_{\Omega}(s)\bigr)_jw^-(s)\, d\sigma_s\,
\]
(cf. \textit{e.g.}, Lanza and Rossi \cite[Thm.~3.1]{LaRo04}.)
Then the Green Identity and classical potential theory imply that
\[
\int_{\partial \Omega}\bigl(\nu_{\Omega}(s)\bigr)_jw^+(s)\,d\sigma_s=\int_{\partial \Omega} s_j \frac{\partial w^+}{\partial \nu_{\Omega}}(s)\,d\sigma_s=\int_{\partial \Omega} s_j \frac{\partial w^-}{\partial \nu_{\Omega}}(s)\,d\sigma_s
\]
(cf. \textit{e.g.}, Lanza and Rossi \cite[Thm.~3.1]{LaRo04}.) As a consequence, since $w^-=\tilde{u}$ and $\int_{\partial \Omega}\bigl(\nu_{\Omega}(s)\bigr)_j\, d\sigma_s=0$, we have
\[
\begin{split}
\int_{\partial \Omega}\bigl(\nu_{\Omega}(s)\bigr)_j\tilde{\theta}(s)\, d\sigma_s&=\int_{\partial \Omega}s_j \frac{\partial \tilde{u}}{\partial \nu_\Omega}(s)\, d\sigma_s-\int_{\partial \Omega} (\nu_\Omega(s))_j \tilde{u}(s)\, d\sigma_s \\
&=\int_{\partial \Omega}s_j \frac{\partial \tilde{u}}{\partial \nu_\Omega}(s)\, d\sigma_s-\int_{\partial \Omega} (\nu_\Omega(s))_j g_0(s)\, d\sigma_s\,.
\end{split}
\]
Accordingly \eqref{eq:repbis} holds and so the proof of (i) is complete. We now consider (ii). Let $R>0$ be such that $(\mathrm{cl}\widetilde{V}\cup \mathrm{cl} \Omega)\subseteq \mathbb{B}_n(0,R)$. By the real analyticity of the restriction operator from $C^{m,\alpha}(\mathrm{cl}\mathbb{B}_n(0,R)\setminus \Omega)$ to $C^{m,\alpha}(\mathrm{cl}\widetilde{V})$, it suffices to prove statement (ii) with $\widetilde{V}$ replaced by $\mathbb{B}_n(0,R)\setminus \mathrm{cl}\Omega$. By taking $\tilde{\epsilon}_2 \in ]0,\epsilon_1]$ small enough, we can assume that 
\[
p+\epsilon \mathrm{cl}\mathbb{B}_n(0,R)\subseteq Q \qquad \forall \epsilon \in ]-\tilde{\epsilon}_2,\tilde{\epsilon}_2[\,.
\] 
If $(\epsilon,g) \in ]0,\tilde{\epsilon}_2[\times \mathcal{U}$, a simple computation based on the Theorem of change of variables in integrals shows that
\begin{multline}
u[\epsilon,g](p+\epsilon t)=-\int_{\partial\Omega}(DS_{n}(t-s) )\nu_{\Omega}(s)\Theta[\epsilon,g](s)\,d\sigma_{s}-\epsilon^{n-1}\int_{\partial\Omega}(DR_{n}^{q}(\epsilon(t-s)) )\nu_{\Omega}(s)\Theta[\epsilon,g](s)\,d\sigma_{s} +\Xi[\epsilon,g] \\\qquad \forall t \in \mathrm{cl}\mathbb{B}_n(0,R) \setminus \mathrm{cl}\Omega\, .\nonumber
\end{multline}
Now we set
\[
\widetilde{G}_{R_n^q}[\epsilon,g](t)\equiv-\epsilon^{n-1}\int_{\partial\Omega}(DR_{n}^{q}(\epsilon(t-s)) )\nu_{\Omega}(s)\Theta[\epsilon,g](s)\,d\sigma_{s} \qquad \forall t \in \mathrm{cl}\mathbb{B}_n(0,R) \setminus \Omega\, ,
\]
for all $(\epsilon,g) \in ]-\tilde{\epsilon}_2,\tilde{\epsilon}_2[\times \mathcal{U}$. Then we note that
\[
\epsilon\mathrm{cl}\mathbb{B}_n(0,R)-\epsilon \partial \Omega \subseteq (\mathbb{R}^n \setminus q\mathbb{Z}^n)\cup \{0\} \qquad \forall \epsilon \in ]-\tilde{\epsilon}_2,\tilde{\epsilon}_2[\,.
\]
Accordingly, by arguing as in the proof of Proposition \ref{prop:Lambda}, we can conclude that $\widetilde{G}_{R_n^q}[\cdot,\cdot]$ is real analytic from $]-\tilde{\epsilon}_2,\tilde{\epsilon}_2[\times \mathcal{U}$ to $C^{m,\alpha}(\mathrm{cl}\mathbb{B}_n(0,R) \setminus \Omega)$. Moreover, $\widetilde{G}_{R_n^q}[0,g_0](\cdot)=0$ in $\mathrm{cl}\mathbb{B}_n(0,R) \setminus \Omega$. By classical results of potential theory and by the real analyticity of $\Theta[\cdot,\cdot]$, there exists a real analytic map $\widetilde{G}_{S_n}[\cdot,\cdot]$ from $]-\tilde{\epsilon}_2,\tilde{\epsilon}_2[\times \mathcal{U}$ to $C^{m,\alpha}(\mathrm{cl}\mathbb{B}_n(0,R) \setminus \Omega)$, such that
\[
\widetilde{G}_{S_n}[\epsilon,g](t)=-\int_{\partial\Omega}(DS_{n}(t-s) )\nu_{\Omega}(s)\Theta[\epsilon,g](s)\,d\sigma_{s}\qquad \forall t \in \mathrm{cl}\mathbb{B}_n(0,R) \setminus \mathrm{cl}\Omega\, ,
\]
for all $(\epsilon,g) \in ]-\tilde{\epsilon}_2,\tilde{\epsilon}_2[\times \mathcal{U}$ (cf.~\textit{e.g.}, Miranda \cite{Mi65}, Dalla Riva and Lanza \cite[Theorem 3.1]{DaLa10a}, Lanza and Rossi \cite[Theorem 3.1]{LaRo04}.) In particular,
\[
\widetilde{G}_{S_n}[0,g_0](t)=\tilde{u}(t)\qquad \forall t \in \mathrm{cl}\mathbb{B}_n(0,R) \setminus \Omega\, .
\]
Then we set
\[
\widetilde{U}[\epsilon,g](t)\equiv \widetilde{G}_{S_n}[\epsilon,g](t)+\widetilde{G}_{R_n^q}[\epsilon,g](t) \qquad \forall t \in \mathrm{cl}\mathbb{B}_n(0,R) \setminus \Omega\, ,
\]
for all $(\epsilon,g) \in ]-\tilde{\epsilon}_2,\tilde{\epsilon}_2[\times \mathcal{U}$. As a consequence, $\widetilde{U}$ is a real analytic map from $]-\tilde{\epsilon}_2,\tilde{\epsilon}_2[\times \mathcal{U}$ to $C^{m,\alpha}(\mathrm{cl}\mathbb{B}_n(0,R)\setminus \Omega)$, such that (jj') holds with $\widetilde{V}$ replaced by $\mathbb{B}_n(0,R)\setminus \mathrm{cl}\Omega$. Thus the proof is complete.\hfill $\Box$
\vspace{\baselineskip}

\begin{rem}
Here we observe that Theorem \ref{thm:rep} (i) concerns what can be called the ``macroscopic'' behaviour of the solution, while Theorem \ref{thm:rep} (ii) describes the ``microscopic'' behaviour. Indeed, in Theorem \ref{thm:rep} (i), we consider a bounded open subset $V$ such that $\mathrm{cl}V \subseteq \mathbb{R}^n \setminus (p+q\mathbb{Z}^n)$, \textit{i.e.} such that its closure $\mathrm{cl}V$ does not intersect the set of points in which the holes degenerate when $\epsilon$ goes to $0$. Then, for $\epsilon$ small enough, $\mathrm{cl}V$ is ``far'' from the union of the holes $\bigcup_{z \in \mathbb{Z}^n}(qz+\Omega_\epsilon)$, and we prove a real analytic continuation result for the restriction of the solution to $\mathrm{cl}V$. Instead, in Theorem \ref{thm:rep} (ii) we take a bounded open subset $\widetilde{V}$ of $\mathbb{R}^n\setminus \mathrm{cl}\Omega$ and we consider the behaviour of the restriction of the solution to the set $p+\epsilon \mathrm{cl}\widetilde{V}$. We note that the set $p+\epsilon \mathrm{cl}\widetilde{V}$ gets, in a sense, closer to the hole $\Omega_\epsilon$ as $\epsilon$ goes to $0$, and that it degenerates into the set $\{p\}$ for $\epsilon=0$. Therefore, in Theorem \ref{thm:rep} (ii) we characterize the behaviour of the solution in proximity of the hole $\Omega_\epsilon$ in the fundamental cell $Q$.
\end{rem}

We now turn to consider the energy integral of the solution on a perforated cell, and we prove the following.

\begin{thm}\label{thm:en}
Let $\alpha \in ]0,1[$. Let $m \in \mathbb{N}\setminus \{0\}$. Let \eqref{ass}-\eqref{assg} hold. Let $\epsilon_1$, $\mathcal{U}$ be as in Proposition \ref{prop:Lambda} (iv). Then there exist $\epsilon_3 \in ]0,\epsilon_1]$ and a real analytic map $G$ from $]-\epsilon_3,\epsilon_3[\times \mathcal{U}$ to $\mathbb{R}$, such that
\begin{equation}\label{eq:en1}
\int_{Q \setminus \mathrm{cl}\Omega_\epsilon}|D_x u[\epsilon,g](x)|^2\, dx=\epsilon^{n-2}G[\epsilon,g]\, ,
\end{equation}
for all $(\epsilon,g) \in ]0,\epsilon_3[\times \mathcal{U}$. Moreover,
\begin{equation}\label{eq:en2}
G[0,g_0]=\int_{\mathbb{R}^n \setminus \mathrm{cl}\Omega}|D\tilde{u}(t)|^2\, dt\,,
\end{equation}
where $\tilde{u}$ is as in Lemma \ref{lem:limeq}.
\end{thm}
{\bf Proof.} Let $(\epsilon,g) \in ]0,\epsilon_1[\times \mathcal{U}$. By the Green Formula and by the periodicity of $u[\epsilon,g](\cdot)$, we have
\begin{equation}\label{eq:en3}
\begin{split}
\int_{Q \setminus \mathrm{cl}\Omega_\epsilon}&|D_x u[\epsilon,g](x)|^2\, dx=\int_{\partial Q}D_x u[\epsilon,g](x)\nu_{Q}(x)u[\epsilon,g](x)\, d\sigma_x-\int_{\partial \Omega_\epsilon}D_x u[\epsilon,g](x)\nu_{\Omega_\epsilon}(x)u[\epsilon,g](x)\, d\sigma_x\\
&=-\epsilon^{n-1}\int_{\partial \Omega}D_x u[\epsilon,g](p+\epsilon t)\nu_{\Omega}(t)g(t)\, d\sigma_t=-\epsilon^{n-2}\int_{\partial \Omega}D \bigl(u[\epsilon,g]\circ(p+\epsilon \mathrm{id}_n)\bigr)(t)\nu_{\Omega}(t)g(t)\, d\sigma_t\,.
\end{split}
\end{equation}
Let $R>0$ be such that $\mathrm{cl}\Omega \subseteq \mathbb{B}_n(0,R)$. By Theorem \ref{thm:rep} (ii), there exist $\epsilon_3 \in ]0,\epsilon_1]$ and a real analytic map $\widetilde{G}[\cdot,\cdot]$ from $]-\epsilon_3,\epsilon_3[\times \mathcal{U}$ to $C^{m,\alpha}(\mathrm{cl}\mathbb{B}_n(0,R)\setminus \Omega)$, such that
\[
p+\epsilon \mathrm{cl}(\mathbb{B}_n(0,R)\setminus \mathrm{cl}\Omega) \subseteq Q\setminus \Omega_\epsilon \qquad \forall \epsilon \in ]-\epsilon_3,\epsilon_3[\setminus \{0\}\, ,
\]
and that
\[
\widetilde{G}[\epsilon,g](t)=u[\epsilon,g]\circ (p+\epsilon \mathrm{id}_n)(t) \qquad\forall t \in \mathrm{cl}\mathbb{B}_n(0,R)\setminus \Omega \quad \forall (\epsilon,g) \in ]0,\epsilon_3[\times \mathcal{U}\, ,
\]
and that
\[
\widetilde{G}[0,g_0](t)=\tilde{u}(t)+\tilde{\xi} \qquad \forall t \in \mathrm{cl}\mathbb{B}_n(0,R)\setminus \Omega\,,
\]
where $\tilde{u}$, $\tilde{\xi}$ are as in Lemma \ref{lem:limeq}. By equality \eqref{eq:en3} we have
\[
\int_{Q \setminus \mathrm{cl}\Omega_\epsilon}|D_x u[\epsilon,g](x)|^2\, dx=-\epsilon^{n-2}\int_{\partial \Omega}D_t \widetilde{G}[\epsilon,g](t)\nu_{\Omega}(t)g(t)\, d\sigma_t\,,
\]
for all $(\epsilon,g) \in ]0,\epsilon_3[\times \mathcal{U}$. Thus it is natural to set
\[
G[\epsilon,g]\equiv -\int_{\partial \Omega}D_t \widetilde{G}[\epsilon,g](t)\nu_{\Omega}(t)g(t)\, d\sigma_t\,,
\]
for all $(\epsilon,g) \in ]-\epsilon_3,\epsilon_3[\times \mathcal{U}$. Then by continuity of the partial derivatives from $C^{m,\alpha}(\mathrm{cl}\mathbb{B}_n(0,R) \setminus \Omega)$ to $C^{m-1,\alpha}(\mathrm{cl}\mathbb{B}_n(0,R) \setminus \Omega)$, and by continuity of the trace operator on $\partial \Omega$ from $C^{m-1,\alpha}(\mathrm{cl}\mathbb{B}_n(0,R) \setminus \Omega)$ to $C^{m-1,\alpha}(\partial \Omega)$, and by the continuity of the pointwise product in Schauder spaces, and by standard calculus in Banach spaces, we conclude that $G[\cdot,\cdot]$ is a real analytic map from $]-\epsilon_3,\epsilon_3[\times \mathcal{U}$ to $\mathbb{R}$ and that equality \eqref{eq:en1} holds. Finally, we note that
\[G[0,g_0]=-\int_{\partial \Omega}D\tilde{u}(t)\nu_{\Omega}(t)g_0(t)\, d\sigma_t\, .
\]
By classical potential theory and the Divergence Theorem, we have
\begin{equation}\label{eq:en4}
\int_{\partial \Omega}D\tilde{u}(t)\nu_{\Omega}(t)\, d\sigma_t\,=0\,. 
\end{equation}
Then, by the decay properties at infinity of $\tilde{u}$ and of its radial derivative and by \eqref{eq:en4}, we have
\[
-\int_{\partial \Omega}D\tilde{u}(t)\nu_{\Omega}(t)g_0(t)\, d\sigma_t=-\int_{\partial \Omega}D\tilde{u}(t)\nu_{\Omega}(t)\bigl(g_0(t)-\tilde{\xi}\bigr)\, d\sigma_t=\int_{\mathbb{R}^n \setminus \mathrm{cl}\Omega}|D\tilde{u}(t)|^2\, dt\, 
\]
(cf.~\textit{e.g.}, Folland \cite[p.~118]{Fo95}.) As a consequence, equality \eqref{eq:en2} follows and the proof is complete. \hfill $\Box$
\vspace{\baselineskip}

\begin{rem}
In Theorem \ref{thm:en}, we have shown a real analytic continuation result for the energy integral of the solution on the perforated fundamental cell $Q \setminus \mathrm{cl}\Omega_\epsilon$, which degenerates into the set $Q \setminus \{p\}$ for $\epsilon=0$. We note that the energy integral $\int_{Q \setminus \mathrm{cl}\Omega_\epsilon}|D_x u[\epsilon,g](x)|^2\, dx$ tends to $0$ as $(\epsilon,g)$ goes to $(0,g_0)$ if $n\geq 3$, while in general this is not true if $n=2$. Moreover, since the map from $]-\epsilon_3,\epsilon_3[$ to $\mathbb{R}$ which takes $\epsilon$ to $\epsilon^{n-2}G[\epsilon,g_0]$ is real analytic, Theorem \ref{thm:en} implies the existence of $\epsilon_3^{\#}\in ]0,\epsilon_3]$ and of a sequence of real numbers $\{a_j\}_{j=0}^{\infty}$ such that 
\[
\int_{Q \setminus \mathrm{cl}\Omega_\epsilon}|D_x u[\epsilon,g_0](x)|^2\, dx=\sum_{j=0}^{\infty}a_j \epsilon^j \qquad \forall \epsilon \in ]0,\epsilon_3^\#[\, ,
\]
where the series converges absolutely in $]-\epsilon_3^\#,\epsilon_3^\#[$. Clearly, analogous considerations for the ``macroscopic'' and ``microscopic'' behaviour of the solution can be derived from the results of Theorem \ref{thm:rep}.
\end{rem}

\appendix
\section{Appendix}\label{app}

In this Appendix, we collect some results exploited in the article.

We have the following known consequence of the Maximum Principle.
\begin{prop}\label{prop:maxprin}
Let $\mathbb{I}$ be a bounded connected open subset of $\mathbb{R}^n$ such that  $\mathbb{R}^n \setminus \mathrm{cl}\mathbb{I}$ is connected and that $\mathrm{cl}\mathbb{I} \subseteq Q$. Let $u \in C^0(\mathrm{cl}\mathbb{S}[\mathbb{I}]^{-}) \cap C^2(\mathbb{S}[\mathbb{I}]^{-})$ be such that
\[
u(x+qe_i)=u(x) \qquad \forall x \in \mathrm{cl}\mathbb{S}[\mathbb{I}]^{-}, \quad \forall i \in \{1,\dots,n\}\, ,
\]
and that
\[
\Delta u(x)= 0 \qquad \forall x \in \mathbb{S}[\mathbb{I}]^{-}\, .
\]
Then the following statements hold.
\begin{enumerate}
\item[(i)]
If there exists a point $x_0 \in \mathbb{S}[\mathbb{I}]^{-}$ such that $u(x_0)=\max_{\mathrm{cl} \mathbb{S}[\mathbb{I}]^{-}} u$, then $u$ is constant within $\mathbb{S}[\mathbb{I}]^{-}$.
\item[(ii)]
If there exists a point $x_0 \in \mathbb{S}[\mathbb{I}]^{-}$ such that $u(x_0)=\min_{\mathrm{cl}\mathbb{S}[\mathbb{I}]^{-}} u$, then $u$ is constant within $\mathbb{S}[\mathbb{I}]^{-}$.
\item[(iii)]
\[
\max_{\mathrm{cl}\mathbb{S}[\mathbb{I}]^{-}}u = \max_{\partial \mathbb{I}} u\, ,
\qquad
\min_{\mathrm{cl}\mathbb{S}[\mathbb{I}]^{-}}u = \min_{\partial \mathbb{I}} u\, .
\]
\end{enumerate}
\end{prop}
{\bf Proof.} Clearly, statement (iii) is a straightforward consequence of (i) and (ii). Furthermore, statement (ii) follows from statement (i) by replacing $u$ with $-u$. Therefore, it suffices to prove (i). Let $u$ and $x_0$ be as in the hypotheses. By periodicity of $u$, $\sup_{x \in \mathbb{S}[\mathbb{I}]^-}u(x)<+\infty$. Then by the Maximum Principle, $u$ must be constant in $\mathbb{S}[\mathbb{I}]^-$ (cf. \textit{e.g.}, Folland \cite[Theorem 2.13, p. 72]{Fo95}.)
 \hfill $\Box$
\vspace{\baselineskip}

We now introduce the following Proposition on nonlinear integral operators (see \cite{LaMu10b}.)

\begin{prop}\label{prop:inop}
Let $n, s\in{\mathbb{N}}$,   $1\leq s< n$. Let ${\mathbb{M}}$ be a compact manifold of class $C^{1}$ imbedded into ${\mathbb{R}}^{n}$ and of dimension $s$. Let ${\mathcal{K}}$ be a Banach space. Let $\mathcal{W}$ be an open subset of $\mathbb{R}^n \times \mathbb{R}^n \times \mathcal{K}$. Let $G$ be a real analytic map from $\mathcal{W}$ to ${\mathbb{R}}$. Then the following statements hold.
\begin{enumerate}
\item[(i)] Let $r \in{\mathbb{N}}$. Let $\Omega$ be a bounded open subset of ${\mathbb{R}}^{n}$. Let 
\[
\mathcal{F}\equiv \Bigl\{(\phi,z) \in C^0(\mathbb{M},\mathbb{R}^n)\times \mathcal{K} \colon \mathrm{cl}\Omega \times\phi(\mathbb{M})\times \{z\}\subseteq \mathcal{W}\Bigr\}\,.
\]
Then the map $H_G$ from $\mathcal{F}\times L^{1}({\mathbb{M}})$ to $C^{r}({\mathrm{cl}}\Omega)$ defined by 
\[
H_G[\phi,z,f](x)\equiv\int_{{\mathbb{M}}}G(x,\phi(y),z)f (y)\,d\sigma_{y}\qquad\forall x\in {\mathrm{cl}}\Omega \,,
\]
for all $(\phi,z,f)\in \mathcal{F}\times L^{1}({\mathbb{M}})$ is real analytic.
\item[(ii)] Let $m \in{\mathbb{N}}$. Let $\alpha \in ]0,1]$. Let $\Omega '$ be a bounded connected open subset of $\mathbb{R}^n$ of class $C^1$. Let 
\[
\mathcal{F}^\#\equiv \Bigl\{(\psi,\phi,z) \in C^{m,\alpha}(\mathrm{cl}\Omega',\mathbb{R}^n)\times C^0(\mathbb{M},\mathbb{R}^n)\times \mathcal{K} \colon \psi(\mathrm{cl}\Omega') \times\phi(\mathbb{M})\times \{z\}\subseteq \mathcal{W}\Bigr\}\,.
\]
Let $H_G^\#$ be the map from $\mathcal{F}^\#\times L^{1}({\mathbb{M}})$ to $C^{m,\alpha}({\mathrm{cl}}\Omega')$ defined by 
\[
H_G^\#[\psi,\phi,z,f](t)\equiv\int_{{\mathbb{M}}}G(\psi(t),\phi(y),z)f (y)\,d\sigma_{y}\qquad\forall t \in {\mathrm{cl}}\Omega' \,,
\]
for all $(\psi,\phi,z,f)\in \mathcal{F}^\#\times L^{1}({\mathbb{M}})$. Then $H_G^\#$ is real analytic from $\mathcal{F}^\#\times L^{1}({\mathbb{M}})$ to $C^{m,\alpha}({\mathrm{cl}}\Omega')$.
\end{enumerate}
\end{prop}

Then we have the following result of (periodic) potential theory (see also Lanza \cite[p. 283]{La10}, Kirsch \cite{Ki89}.)
\begin{prop}\label{prop:comphom}
Let  $\alpha\in]0,1[$, $m\in {\mathbb{N}}\setminus\{0\}$. Let $S_{n}^{q}$ be as in Theorem  \ref{psper}. Let $\mathbb{I}$ be a bounded connected open subset of ${\mathbb{R}}^{n}$ of class $C^{m,\alpha}$ such that 
${\mathbb{R}}^{n}\setminus{\mathrm{cl}}\mathbb{I}$ is connected and that ${\mathrm{cl}}\mathbb{I}\subseteq Q$. Let $S_{n}^{q}$ be as in Theorem  \ref{psper}. Then the following statements hold.
\begin{enumerate}
\item[(i)] The map $w_q[\partial \mathbb{I},\cdot]_{|\partial \mathbb{I}}$ is compact from $C^{m,\alpha}(\partial \mathbb{I})$ to itself.
\item[(ii)] Let $\widetilde{M}[\cdot]$ be the map from $C^{m,\alpha}(\partial \mathbb{I})$ to itself, defined by
\[
\widetilde{M}[\mu](t)\equiv -\frac{1}{2}\mu(t)+w_{q}[\partial \mathbb{I},\mu](t)\qquad \forall t \in \partial \mathbb{I}\, ,
\]
for all $\mu \in C^{m,\alpha}(\partial \mathbb{I})$. Then $\widetilde{M}[\cdot]$ is a linear homeomorphism from $C^{m,\alpha}(\partial \mathbb{I})$ onto itself. Moreover,
\begin{equation}\label{eq:comphom}
\Bigl\{\widetilde{M}[\lambda]\colon \lambda \in \mathbb{R} \Bigr\}=\mathbb{R}\,,
\end{equation}
where we identify the constant functions with the constants themselves.
\end{enumerate}
\end{prop}
 {\bf Proof.} We start by proving (i). Let $R_{n}^{q}$ be as in Theorem  \ref{psper}. We set
 \[
 w[\partial\mathbb{I},\mu](t)\equiv
-\int_{\partial\mathbb{I}}(DS_{n}(t-s) )\nu_{\mathbb{I}}(s)\mu(s)\,d\sigma_{s}
\qquad\forall t\in \partial \mathbb{I}\,,
 \]
 for all $\mu \in C^{m,\alpha}(\partial \mathbb{I})$. By classical potential theory and by the compactness of the imbedding of $C^{m,\alpha}(\partial \mathbb{I})$ into $C^{m,\beta}(\partial \mathbb{I})$ for $\beta \in ]0,\alpha[$, we conclude that the operator $w[\partial \mathbb{I},\cdot]$ from $C^{m,\alpha}(\partial \mathbb{I})$ to itself is compact. Indeed, case $m=1$ has been proved by Schauder \cite[Hilfsatz XI, p.~618]{Sc31}, and case $m>1$ follows by taking the tangential derivatives of $w[\partial \mathbb{I},\cdot]$ on $\partial \mathbb{I}$ and by arguing by induction on $m$. We also set
\[
 w_{R_{n}^q}[\partial\mathbb{I},\mu](t)\equiv
-\int_{\partial\mathbb{I}}(DR_{n}^{q}(t-s) )\nu_{\mathbb{I}}(s)\mu(s)\,d\sigma_{s}
\qquad\forall t\in \partial \mathbb{I}\,,
 \]
 for all $\mu \in C^{m,\alpha}(\partial \mathbb{I})$. Clearly, $w_q[\partial \mathbb{I},\mu]=w[\partial \mathbb{I},\mu]+w_{R_n^q}[\partial \mathbb{I},\mu]$ on $\partial \mathbb{I}$, for all $\mu \in C^{m,\alpha}(\partial \mathbb{I})$. For each $j \in \{1,\dots,n\}$, we set
 \[
 N_{R_{n}^{q},j}[f](t)\equiv
-\int_{\partial\mathbb{I}}(D_{x_j}R_{n}^{q}(t-s) )f(s)\,d\sigma_{s}
\qquad\forall t\in \mathrm{cl} \mathbb{I}\,,
 \]
 for all $f \in L^1(\partial \mathbb{I})$. By Proposition \ref{prop:inop} (i), $N_{R_{n}^{q},j}[\cdot]$ is linear and continuous from $L^1(\partial \mathbb{I})$ to $C^{m+1}(\mathrm{cl}\mathbb{I})$. Moreover, by the compactness of the imbedding of $C^{m+1}(\mathrm{cl} \mathbb{I})$ into $C^{m,\alpha}(\mathrm{cl} \mathbb{I})$, $N_{R_{n}^{q},j}[\cdot]$ is compact from $L^1(\partial \mathbb{I})$ to $C^{m,\alpha}(\mathrm{cl}\mathbb{I})$ (cf.~\textit{e.g.}, Lanza and Rossi \cite[Lemma 2.1]{LaRo04}.) Then by the continuity of the map from $C^{m,\alpha}(\partial \mathbb{I})$ to $L^1(\partial \mathbb{I})$ which takes $\mu$ to $(\nu_{\mathbb{I}})_j\mu$, and by the continuity of the trace operator from $C^{m,\alpha}(\mathrm{cl} \mathbb{I})$ to $C^{m,\alpha}(\partial \mathbb{I})$, we immediately deduce the compactness of $w_{R_{n}^{q}}[\partial \mathbb{I},\cdot]$ from $C^{m,\alpha}(\partial \mathbb{I})$ to $C^{m,\alpha}(\partial \mathbb{I})$, and, as a consequence, of $w_q[\partial \mathbb{I},\cdot]_{|\partial \mathbb{I}}$. Hence the proof of (i) is complete. We now turn to the proof of (ii). By the Open Mapping Theorem, it suffices to prove that $\widetilde{M}[\cdot]$ is a bijection. By (i) and by the Fredholm Theory, it suffices to show that $\widetilde{M}[\cdot]$ is injective. So let $\mu \in C^{m,\alpha}(\partial \mathbb{I})$ be such that
\[
-\frac{1}{2}\mu+w_q[\partial \mathbb{I},\mu]=0 \qquad \text{on $\partial \mathbb{I}$}\,.
\]
 By Theorem \ref{dperpot}, $w_q^-[\partial \mathbb{I},\mu]$ is a solution of the following problem
 \[
  \left \lbrace 
 \begin{array}{ll}
 \Delta u (x)= 0 & \textrm{$\forall x \in {\mathbb{S}} [\mathbb{I}]^{-}$}\,, \\
u(x+qe_i) =u(x) &  \textrm{$\forall x \in \mathrm{cl} {\mathbb{S}} [\mathbb{I}]^{-}, \quad \forall i \in \{1,\dots,n\}$}\,, \\
u(x)=0& \textrm{$\forall x \in \partial \mathbb{I}$}\,.
 \end{array}
 \right.
\]
 As a consequence, by Proposition \ref{prop:maxprin}, $w_q^-[\partial \mathbb{I},\mu]= 0$ in $\mathrm{cl}\mathbb{S}[\mathbb{I}]^{-}$. In particular, 
 \[
 \frac{\partial}{\partial \nu_{\mathbb{I}}}w_q^-[\partial \mathbb{I},\mu]=0 \qquad \text{on $\partial \mathbb{I}$}\,.
 \]
 Then, by formula \eqref{dperpot2b},
 \[
 \frac{\partial}{\partial \nu_{\mathbb{I}}}w_q^+[\partial \mathbb{I},\mu]=0 \qquad \text{on $\partial \mathbb{I}$}\,.
 \]
 Accordingly, by Theorem \ref{dperpot}, $w_q^+[\partial \mathbb{I},\mu]_{|\mathrm{cl}\mathbb{I}}\in C^{m,\alpha}(\mathrm{cl}\mathbb{I})$ is a solution of the following problem
  \[
  \left \lbrace 
 \begin{array}{ll}
 \Delta u (x)= 0 & \textrm{$\forall x \in \mathbb{I}$}\,, \\
\frac{\partial}{\partial \nu_{\mathbb{I}}}u(x)=0& \textrm{$\forall x \in \partial \mathbb{I}$}\,.
 \end{array}
 \right.
\]
 As a consequence, there exists a constant $c\in \mathbb{R}$ such that $w_q^+[\partial \mathbb{I},\mu]=c$ on $\mathrm{cl}\mathbb{S}[\mathbb{I}]$. By formula \eqref{dperpot2a},
 \[
 \mu=w_q^+[\partial \mathbb{I},\mu]-w_q^-[\partial \mathbb{I},\mu]=c \qquad \text{on $\partial \mathbb{I}$}\,.
 \]
Therefore, by formula \eqref{dperpot3c},
\[
\widetilde{M}[\mu]=w_q^-[\partial \mathbb{I},c]=-c\frac{|\mathbb{I}|}{|Q|} \qquad \text{on $\partial \mathbb{I}$}\,,
\]
and so $c=0$. Hence, $\mu=0$. Finally, equality \eqref{eq:comphom} follows immediately from \eqref{dperpot3a}. Thus the proof is complete.
 \hfill $\Box$
\vspace{\baselineskip}

Finally, we have the following well known result of classical potential theory.

\begin{prop}\label{prop:classpot}
Let  $\alpha\in]0,1[$, $m\in {\mathbb{N}}\setminus\{0\}$. Let $\Omega$ be a bounded connected open subset of ${\mathbb{R}}^{n}$ of class $C^{m,\alpha}$. Let $\widetilde{N}[\cdot,\cdot]$ be the map from $C^{m,\alpha}(\partial \Omega)_0 \times \mathbb{R}$ to $C^{m,\alpha}(\partial \Omega)$, defined by
\[
\widetilde{N}[\mu,\xi](x)\equiv -\frac{1}{2}\mu(x)-\int_{\partial\Omega}(DS_{n}(x-y) )\nu_{\Omega}(y)\mu(y)\,d\sigma_{y}+\xi \qquad \forall x \in \partial \Omega\, ,
\]
for all $(\mu,\xi) \in C^{m,\alpha}(\partial \Omega)_0\times \mathbb{R}$. Then $\widetilde{N}[\cdot,\cdot]$ is a linear homeomorphism from $C^{m,\alpha}(\partial \Omega)_0\times \mathbb{R}$ onto $C^{m,\alpha}(\partial \Omega)$. 
\end{prop}
 {\bf Proof.} Clearly, $\widetilde{N}$ is linear and continuous (cf.~\textit{e.g.}, Miranda \cite{Mi65}, Dalla Riva and Lanza \cite[Theorem 3.1]{DaLa10a}, Lanza and Rossi \cite[Theorem 3.1]{LaRo04}.) By the Open Mapping Theorem, it suffices to show that it is a bijection. By well known results of classical potential theory, we have
\[
C^{m,\alpha}(\partial \Omega)=\Bigl\{-\frac{1}{2}\mu(\cdot)-\int_{\partial\Omega}(DS_{n}(\cdot-y) )\nu_{\Omega}(y)\mu(y)\,d\sigma_{y} \colon \mu \in C^{m,\alpha}(\partial \Omega)\Bigr\}\oplus <\chi_{\partial \Omega}>\, ,
\]
where $\chi_{\partial \Omega}$ denotes the characteristic of $\partial \Omega$ (cf.~\textit{e.g.}, Folland \cite[Ch.~3]{Fo95} and Lanza \cite[Appendix A]{La07a}.) On the other hand, as is well known, for each $\psi$ in the set
\[
\Bigl\{-\frac{1}{2}\mu(\cdot)-\int_{\partial\Omega}(DS_{n}(\cdot-y) )\nu_{\Omega}(y)\mu(y)\,d\sigma_{y} \colon \mu \in C^{m,\alpha}(\partial \Omega)\Bigr\}\, ,
\]
there exists a unique $\mu$ in $C^{m,\alpha}(\partial \Omega)$ such that
\[
  \left \lbrace 
 \begin{array}{ll}
\psi(x)=-\frac{1}{2}\mu(x)-\int_{\partial\Omega}(DS_{n}(x-y) )\nu_{\Omega}(y)\mu(y)\,d\sigma_{y} \qquad \forall x \in \partial \Omega\, ,\\
\int_{\partial \Omega}\mu \, d\sigma=0 & 
 \end{array}
 \right.
\]
(cf.~\textit{e.g.}, Folland \cite[Ch.~3]{Fo95} and Lanza \cite[Appendix A]{La07a}.) As a consequence, for each $\phi \in C^{m,\alpha}(\partial \Omega)$, there exists a unique pair $(\mu,\xi)$ in $C^{m,\alpha}(\partial \Omega)_0 \times \mathbb{R}$, such that
\[
\phi(x)=-\frac{1}{2}\mu(x)-\int_{\partial\Omega}(DS_{n}(x-y) )\nu_{\Omega}(y)\mu(y)\,d\sigma_{y}+\xi \qquad \forall x \in \partial \Omega\, ,
\]
and so $\widetilde{N}$ is bijective. Thus the proof is complete.
 \hfill $\Box$
\vspace{\baselineskip}

\section*{Acknowledgements} This paper generalizes a part of the work performed by the author in his ``Laurea Specialistica'' Thesis \cite{Mu08} under the guidance of  Prof. M. Lanza de Cristoforis. The author wishes to thank Prof. M. Lanza de Cristoforis for his constant help during the preparation of this paper. The results presented here have been announced in \cite{Mu10}. The author acknowledges the support of the research project ``Un approccio funzionale analitico per problemi di omogeneizzazione in domini a perforazione periodica'' of the University of Padova, Italy.

\end{document}